\documentclass[11pt]{amsart}
\usepackage[T1]{fontenc}
\usepackage{lmodern,amssymb,mathtools,booktabs,array,microtype,etoolbox}
\usepackage{tikz}
\usetikzlibrary{arrows.meta,decorations.markings}
\tikzset{midarrow/.style={postaction={decorate},decoration={markings,mark=at position 0.55 with {\arrow{Stealth[length=2mm]}}}}}
\usepackage[margin=1.1in]{geometry}
\usepackage[colorlinks=true,linkcolor=blue,citecolor=blue,urlcolor=blue]{hyperref}
\makeatletter
\patchcmd{\section}{\scshape}{\bfseries}{}{}
\makeatother
\allowdisplaybreaks
\newtheorem{theorem}{Theorem}[section]
\newtheorem{proposition}[theorem]{Proposition}
\newtheorem{lemma}[theorem]{Lemma}
\newtheorem{corollary}[theorem]{Corollary}
\newtheorem{maintheorem}{Theorem}

\theoremstyle{remark}
\newtheorem{remark}[theorem]{Remark}
\newcommand{\HH}{\mathbb H}
\newcommand{\CC}{\mathbb C}
\newcommand{\RR}{\mathbb R}
\newcommand{\ZZ}{\mathbb Z}
\DeclareMathOperator{\Rea}{Re}
\DeclareMathOperator{\Ima}{Im}
\DeclareMathOperator{\sgn}{sgn}

\title[Radial reduction and spherical splitting]{Radial reduction, spherical splitting and root counts for two-sided quaternionic polynomials}
\author{R. Laniewski}
\date{October 14, 2026}

\begin{document}

\begin{abstract}
We study quaternionic polynomials whose coefficients may occur on both sides of the variable, using their degree in the four real coordinates. Such polynomials can have no roots, or many more isolated roots than their degree. A simple root has an invertible real derivative, and its sign is the sign of that derivative's determinant.

For a radial family, the roots over a non-real target reduce to one real equation. We obtain a sharp root bound and follow a change in the signed root count as roots escape to infinity.

The main results concern spherical root sets for which the polynomial derivative is invertible in the two normal directions. Under a small polynomial perturbation, a sphere splits into simple roots corresponding to the zeros of a tangent vector field, when the derivative of that field is invertible at every zero. We obtain a sharp bound in terms of the perturbation degree and a signed count of two or of minus two.

For even degrees, the construction yields polynomials with a single leading power and lower-degree terms whose simple root counts, divided by the fourth power of the degree, tend to one as the degree grows. An explicit quartic has 132 simple roots, and this count persists under small coefficient perturbations. For gradient perturbations, the roots correspond to equilibrium directions of an energy on the sphere. This correspondence gives energy thresholds for root disappearance and bounds the critical directions of polynomial orientation profiles in diffusion magnetic resonance imaging.
\end{abstract}

\maketitle

\tableofcontents

\section{Introduction}\label{sec:introduction}
A two-sided quaternionic polynomial can omit zero. For $q=x_0+x_1i+x_2j+x_3k$, the multiplication rules give $ji-ij=-2k$ and $ki-ik=2j$, hence $qi-iq=2(x_3j-x_2k)$. Thus $qi-iq-1$ has real part $-1$ everywhere. We therefore examine the transition between root existence and rootlessness, together with the root counts along it. We look for families in which the image, the number of isolated roots, their local degrees and their persistence under perturbation can all be determined.

Throughout, a two-sided polynomial expression is a finite sum of words with quaternionic coefficients interspersed with the variable. We evaluate each word in its specified order. We call its degree in the four real coordinates the real coordinate degree. Degree refers to this real coordinate degree unless Brouwer degree is explicitly named. We write $d$ for the power exponent, $D$ for a perturbation degree, and $N$ for the total real coordinate degree. These maps differ from polynomials in a formally central indeterminate and from the one-sided polynomials appearing in quaternionic root-existence theorems.

Every non-constant one-sided quaternionic polynomial has a root \cite{niven}. Eilenberg and Niven also allow a single highest-degree word with non-zero coefficients, together with lower-degree terms \cite{en}. The companion-matrix approach to unilateral polynomials is developed in \cite{deleo}. Janovsk\'{a} and Opfer study the classification and numerical computation of zeros of two-sided polynomials~\cite{janovskaopfer}. The exact counts and local degrees below complement this classification for radial and imaginary perturbations. General two-sided expressions admit further possibilities, and their representations and polynomial identities require separate attention \cite{mierzejewski}.

Evaluation of two-sided words gives the full space of real polynomial maps $\RR^4\to\RR^4$, as shown in~\eqref{eq:coefficient-space}. General persistence and B\'ezout estimates in that space are standard real and complex polynomial arguments. The quaternionic structure considered here comes from retaining the power map and perturbing it by a real-valued term, or by an imaginary term that depends only on the imaginary part of the variable.

For a degree-$n$ polynomial in a central indeterminate over a division ring, the Gordon--Motzkin theorem bounds the number of conjugacy classes containing roots by $n$ \cite{gordonmotzkin}. In the quaternionic one-sided case, isolated roots and spherical zero sets are studied by Pogorui and Shapiro \cite{pogoruishapiro}. The radial two-sided expressions below belong to a broader class, so these classical bounds do not directly govern their fibers.
Our point of departure is the class
\begin{equation}\label{eq:class}
P_{d,\phi}(q)=q^d+\phi(|q|^d),\qquad d=2m\geq2,\qquad \phi\in\RR[t].
\end{equation}
The even exponent ensures polynomiality in real coordinates. A constant $\phi$ gives a real translation of the power map $Q_d(q)=q^d$, of real coordinate degree and Brouwer degree $d$, while the bound in Theorem~\ref{thm:count} is stated for $\deg\phi\geq1$. The perturbation depends only on the norm of $q^d$, so the map factors through the quaternionic power map as $P_{d,\phi}=T_\phi\circ Q_d$, with $T_\phi(w)=w+\phi(|w|)$. Both factors respect quaternionic conjugation, and $T_\phi$ preserves the imaginary part. In this family a scalar displacement changes the image and the degree, while the power map keeps its explicit fibers. Centrality localizes non-real fibers to a complex slice. Dependence on the norm then reduces their remaining two real coordinates to one radius, producing the bound $2N$. Section~\ref{sec:central} separates these two reductions, and it shows that a radius measured by $|q|^2$ attains the larger planar bound.

\subsection{Main results}\label{subsec:principal}
A root of a quaternionic polynomial map is called simple when its real derivative is invertible. Its local degree is then the sign of that derivative's determinant, using the orientation $(1,i,j,k)$ on $\HH$. We state the three main results in a short form, and the sections that prove them give the complete statements.

Write a non-real target as $b=a+v$, with $a\in\RR$, $v\in\Ima\HH$, and $\rho=|v|>0$. Set $I=v/\rho$ and $\CC_I=\RR+I\RR$. Define
\begin{equation}\label{eq:radial}
R_b(t)=(a-\phi(t))^2+\rho^2-t^2.
\end{equation}
The radial family is governed by this single real polynomial.

\begin{maintheorem}[The radial family]\label{thm:main-radial}
Let $d=2m\geq2$, let $\phi$ be a real polynomial of degree $r\geq1$, and let $b$ be a non-real target. The roots of $P_{d,\phi}(q)=b$ lie in $\CC_I$, and exactly $d$ of them correspond to each positive root $t$ of $R_b$. They are simple exactly when $R_b'(t)\neq0$, with local degree $-\sgn R_b'(t)$. The fiber has therefore at most $2dr$ points, and for every pair $(d,r)$ this bound is attained with all roots simple.
\end{maintheorem}

Theorem~\ref{thm:count} gives the complete statement, and Section~\ref{sec:sharp} proves the sharpness. In terms of the real coordinate degree $N=dr$ of $P_{d,\phi}$, the bound is $2N$.

The member $q^d+\lambda|q|^d$, in which $\phi$ is linear, shows how the Brouwer degree changes within the family. Theorem~\ref{thm:transition}, stated in Section~\ref{sec:transition}, describes it for every $\lambda\geq0$. Below $\lambda=1$ the map is proper and surjective, of Brouwer degree $d$. Above $\lambda=1$ it is proper of Brouwer degree zero, its image is a closed cone, and the non-zero boundary points of this cone are covered by ordinary folds at which roots of opposite local degree meet. At $\lambda=1$ properness is lost, and the degree changes through this loss.

The radial reduction leaves aside the real targets, whose fibers contain conjugacy two-spheres. These spheres are, in our view, the most interesting part of the problem, since the one-sided and two-sided theories separate most visibly there. Under a small perturbation, a sphere with non-degenerate normal derivative breaks into isolated roots, and Theorem~\ref{thm:sphere-splitting} identifies them with the zeros of a tangent vector field on $\mathbb S^2$. The Poincar\'e--Hopf theorem then fixes their signed count at $\pm2$, and Corollary~\ref{cor:sphere-degree} shows that a sufficiently small perturbation always leaves a root near such a sphere, and at least two roots when all of them are simple. The Lyapunov--Schmidt reduction first solves the real and sphere-parallel equations for the two normal coordinates by the implicit function theorem. After substitution and division by the perturbation parameter, a tangent-section equation remains, and its leading field determines the splitting. Figure~\ref{fig:sphere-splitting} illustrates these stages. The quaternionic structure enters through the form of the reduced field, to which we now turn.

Let $\mathbb S^2$ denote the unit sphere of $\Ima\HH$. Write a conjugacy sphere in a real fiber of $P_{d,\phi}$ as $\Sigma=\{x_0+y_0J\mid J\in\mathbb S^2\}$, with $y_0>0$. On $\Sigma$, one has $q^d=\sigma t_0$, where $\sigma=\pm1$ and $t_0=(x_0^2+y_0^2)^{d/2}$. We call $\Sigma$ normally non-degenerate when the derivative in the two normal coordinates, the real part $x$ and the imaginary radius $y$, is invertible onto the real and $J$-parallel components, as in~\eqref{eq:sphere-normal-matrix}. Lemma~\ref{lem:sphere-normal} gives the equivalent condition $\kappa=1+\sigma\phi'(t_0)\neq0$. For a polynomial perturbation $H$, the first reduced field on $\Sigma$ is the tangential part of $\Ima H$ along the sphere. It is a tangent vector field on $\mathbb S^2$, whose precise form is given in Section~\ref{sec:spherical-splitting}, and a zero of this field is called non-degenerate when the derivative of the field is invertible at that point.

\begin{maintheorem}[Splitting of a conjugacy sphere]\label{thm:main-splitting}
Let $\Sigma$ be a normally non-degenerate conjugacy sphere in the fiber of $P_{d,\phi}$ over a real target $a$. Fix a polynomial map $H$ of real coordinate degree at most $D\geq1$ whose first reduced field on $\Sigma$ has only non-degenerate zeros. For this $H$, choose a sufficiently small isolating tubular neighborhood $U$ of $\Sigma$, independent of $\varepsilon$. There is then $\varepsilon_0(H,U)>0$ such that, for $0<|\varepsilon|<\varepsilon_0(H,U)$, the roots of $P_{d,\phi}+\varepsilon H=a$ contained in $U$ are simple and correspond bijectively to the zeros of that field. Their local degrees add up to $\pm2$, and their number is even and at most $2(D^2+D+1)$. For every $D\geq1$ this bound is attained.
\end{maintheorem}

Theorems~\ref{thm:sphere-splitting} and~\ref{thm:sphere-bounded-count} give the complete statements. For a perturbation $h(\Ima q)$, with $h$ a homogeneous polynomial map of $\Ima\HH$, the first reduced field vanishes exactly at the unit vectors $J$ with $h(J)=\mu J$ for a real number $\mu$. These are the tensor eigenvectors used below. Lemma~\ref{lem:eigen-reduction} makes the reduction exact and global in this case, and the root problem separates into an eigenvector problem on $\Ima\HH$ and a planar problem on a complex slice. The upper bound for all polynomial first reduced fields comes from a Chern class computation on a complex quadric, and homogeneous gradients attain it through the real eigenconfigurations constructed by Kozhasov \cite{kozhasov}.

The eigenvector reduction bears on the theorem of Eilenberg and Niven, which guarantees a root for $q^d$ plus terms of lower degree. A radial term of lower degree creates several conjugacy spheres, and a single imaginary perturbation splits them all at once.

\begin{maintheorem}[One leading word]\label{thm:main-one-word}
For every even $d\geq2$ there is a polynomial map of the form $q^d$ plus terms of real coordinate degree less than $d$ with exactly $d(d-1)(d^2-2d+3)$ roots, all simple. For $d=2$, every map $q^2+L(q)$ with $L$ a real affine map of $\HH$ has at most eight simple roots, and this bound is attained.
\end{maintheorem}

The first assertion follows from Theorem~\ref{thm:en-stacked}, with radii chosen as in the proof of Corollary~\ref{cor:en-quartic}, and the second is Proposition~\ref{prop:en-quadratic}. The lower bound has the same leading coefficient as the B\'ezout bound $d^4$. For $d=4$ it equals $132$, and Corollary~\ref{cor:en-quartic} provides an explicit polynomial. Remark~\ref{rem:en-stacked-signs} divides these roots by local degree.

\subsection{Further results}\label{subsec:further-results}
Table~\ref{tab:root-bounds} compares the root bounds. Proposition~\ref{prop:norm-square} attains the central bound for every even $N\geq d$, and the coordinate grid of Remark~\ref{rem:ambient} attains the ambient bound.
\begin{table}[htbp]
\caption{Root bounds under the hypotheses of the cited results. The radial bound counts distinct roots over a non-real target, and the other bounds count simple roots. The spherical bound concerns one isolating tube. Simple central roots persist locally by the implicit function theorem. Stability of the complete count follows under the coercivity hypothesis of Lemma~\ref{lem:stability}. The one-word entry bounds the maximum simple root count $M_d$ from below and, together with the ambient bound, gives $M_d\sim d^4$ as $d\to\infty$ through even integers. Here $N$ is the real coordinate degree and $D$ is the perturbation degree.}\label{tab:root-bounds}
\small
\begin{tabular}{@{}p{0.38\textwidth}p{0.28\textwidth}p{0.27\textwidth}@{}}
\toprule
Class and scope & Bound & Result\\
\midrule
Radial class, non-real target & At most $2N$ & Theorem~\ref{thm:count}\\
Central class, non-real target & At most $dN$ & Proposition~\ref{prop:central}\\
Full coefficient space & At most $N^4$ & Remark~\ref{rem:ambient}\\
Local spherical splitting & At most $2(D^2+D+1)$ & Theorem~\ref{thm:sphere-bounded-count}\\
One leading word, even $d$ & $M_d\geq d(d-1)(d^2-2d+3)$ & Theorem~\ref{thm:main-one-word}\\
Quadratic one-word class, every fiber & At most eight & Proposition~\ref{prop:en-quadratic}\\
\bottomrule
\end{tabular}
\end{table}

The same reduction describes the Brouwer degree of whole families. For a quadratic imaginary perturbation of $q^2+1$, Proposition~\ref{prop:quad-wall} expresses the degree as a sum of indices over the eigenvectors whose eigenvalues remain below a threshold. Theorem~\ref{thm:sphere-xyz} gives in this way a map of real coordinate degree two with fourteen simple roots, whose degree passes from $2$ to $-6$. Both this transition and that of Theorem~\ref{thm:transition} involve escape to infinity, while gradient perturbations of degree $D\geq3$ exhibit fold collisions at finite radius, as shown in Subsection~\ref{subsec:anisotropy-reach}.

When the imaginary perturbation is the gradient of a homogeneous polynomial, its unit eigenvectors are the equilibrium directions of an anisotropy energy on $\mathbb S^2$. Their Morse types determine the tangent indices, which enter the local-degree formula together with the normal determinant. Section~\ref{sec:anisotropy} develops this correspondence. The biaxial, octupolar and cubic magnetocrystalline energies realize the equality case of Theorem~\ref{thm:sphere-bounded-count} for $D=1,2,3$, and the cubic energy gives an explicit polynomial with $132$ simple roots in Corollary~\ref{cor:en-quartic}.

\subsection{Outline}\label{subsec:outline}
Section~\ref{sec:expressions} identifies two-sided expressions with real polynomial maps and computes the derivative of the power map on a complex slice. Section~\ref{sec:reduction} reduces the fiber of the radial family over a non-real target to the positive roots of one real polynomial, and it reads the local degree of each quaternionic root from the sign of the derivative of that polynomial. Section~\ref{sec:sharp} constructs radial terms for which the resulting bound $2dr$ is attained, with an explicit degree-four radial example of eight simple roots in Proposition~\ref{prop:quartic}. The quartic with $132$ simple roots in Corollary~\ref{cor:en-quartic} has leading word $q^4$. Section~\ref{sec:transition} follows the member $q^d+\lambda|q|^d$ through the value $\lambda=1$, where properness is lost and the Brouwer degree falls from $d$ to zero, and it describes the zero fiber at the threshold and the folds that appear beyond it. Section~\ref{sec:central} compares the radial family with real-valued perturbations that are central without being radial, for which the planar B\'ezout bound $dN$ takes the place of $2N$, and it shows that a radial term in $|q|^2$ attains this larger bound. Section~\ref{sec:stability} supplies the standard coercivity argument by which the extremal counts persist in the full coefficient space.

The second half of the paper turns to the real targets. Section~\ref{sec:spherical-splitting} shows how a conjugacy sphere in a real fiber breaks into isolated roots under a small perturbation, identifies these roots with the zeros of a tangent field on $\mathbb S^2$, and gives a global example in which the count is exact. Section~\ref{sec:eigen} specializes to imaginary perturbations that depend on $\Ima q$ alone. The roots then separate into the unit eigenvectors of a tensor and the solutions of a planar equation, the number of eigenvectors is bounded by $2(D^2+D+1)$, and the same reduction yields polynomials with one leading word $q^d$ whose simple roots grow in number like $d^4$. Section~\ref{sec:anisotropy} reads the eigenvectors as the equilibrium directions of an anisotropy energy and works through the biaxial, octupolar and cubic cases. Section~\ref{sec:mri-orientation} outlines a possible application to orientation finding in diffusion MRI.

\section{Polynomial expressions and quaternionic powers}\label{sec:expressions}
\subsection{The radial term is a two-sided polynomial}\label{subsec:radial-expression}
Quaternionic conjugation has the word expression
\begin{equation}\label{eq:conjugate}
\overline q=-\frac12(q+iqi+jqj+kqk).
\end{equation}
Therefore $|q|^{2m}=(q\overline q)^m$ is a two-sided expression of degree $2m$. Formula~\eqref{eq:class} belongs to the declared class.

For clarity about the coefficient space in Proposition~\ref{prop:stable}, write $q=x_0+x_1i+x_2j+x_3k$. The real coordinate functions satisfy
\[
x_0=\Rea q,\qquad x_1=-\Rea(iq),\qquad x_2=-\Rea(jq),\qquad x_3=-\Rea(kq)
\]
with $\Rea q=(q+\overline q)/2$. Thus every real coordinate polynomial is expressible by quaternionic words, and multiplication by $1,i,j,k$ gives every quaternion-valued real polynomial map. Let $\mathcal W_{\leq N}$ be the real vector space of two-sided expressions of word degree at most $N$. Evaluation induces an isomorphism
\begin{equation}\label{eq:coefficient-space}
\mathcal W_{\leq N}/\ker(\operatorname{ev})
\cong\bigl(\RR[x_0,x_1,x_2,x_3]_{\leq N}\bigr)^4.
\end{equation}
Indeed, evaluation sends each word to a real polynomial map of at most the same degree, and the coordinate formulas above prove surjectivity without increasing degree. The kernel consists of expressions that vanish identically as maps. Thus the perturbation neighborhoods in Proposition~\ref{prop:stable} are open in the full space of polynomial maps, with its real coordinate coefficient topology. Formal expressions are identified here precisely when they define the same function.

\subsection{Powers over a non-real target}\label{subsec:powers}
For a unit imaginary quaternion $I$, write $\Ima_I(x+Iy)=y$, so $\Ima_I$ takes real values on $\CC_I$.

\begin{lemma}\label{lem:power}
Let $d\geq1$ be any integer. For $w\in\HH\setminus\RR$, set $I_w=\Ima w/|\Ima w|$. The equation $q^d=w$ has exactly $d$ solutions, all in $\CC_{I_w}=\RR+I_w\RR$. At each of them the power map $Q_d(q)=q^d$ has positive real Jacobian determinant.
\end{lemma}
\begin{proof}
The equation implies $qw=wq$. The centralizer of a non-real quaternion $w$ is $\CC_{I_w}$, so the solutions are the $d$ distinct complex roots of $w$ in that plane. They are all non-real.

Set $I=I_w$ and write a solution as $q=x+Iy$ with $y\neq0$. The derivative preserves the plane $\CC_I$ and its orthogonal complement. On the plane it is complex multiplication by $d q^{d-1}$. On the orthogonal complement it is real scalar multiplication by $\Ima_I(q^d)/y$. To see the transverse factor explicitly, choose a unit imaginary $J$ orthogonal to $I$, so that $I$, $J$ and $IJ$ obey the multiplication rules of $i$, $j$ and $k$, and write $q=s(\cos\theta+I\sin\theta)$, so $y=s\sin\theta$. Since $JI=-IJ$, one has $J(x+Iy)=(x-Iy)J$. The resulting relation $Jq=\overline qJ$ implies $Jq^{d-1-k}=\overline q^{\,d-1-k}J$, and hence
\begin{align*}
DQ_d(q)[J]&=\sum_{k=0}^{d-1}q^kJq^{d-1-k}
=\left(\sum_{k=0}^{d-1}q^k\overline q^{\,d-1-k}\right)J\\
&=\frac{q^d-\overline q^{\,d}}{q-\overline q}J
=s^{d-1}\frac{\sin(d\theta)}{\sin\theta}J
=\frac{\Ima_I(q^d)}{y}J.
\end{align*}
The quotient is computed in $\CC_I$ and is real. Since $(IJ)I=-I(IJ)$, we also have $(IJ)q=\overline q(IJ)$. The same calculation therefore applies to $IJ$, completing the calculation on the orthogonal complement. With $\mu=\Ima_I(q^d)/y$, the transverse determinant is $\mu^2>0$, independently of the sign of $\mu$. The quaternionic product $IJ$ equals $I\times J$, so the change from $(1,i,j,k)$ to $(1,I,J,IJ)$ lies in $\operatorname{SO}(3)\subset\operatorname{SO}(4)$, with the real axis fixed. The non-real target ensures $\mu\neq0$. Consequently
\begin{equation}\label{eq:powerjac}
\det DQ_d(q)=d^2|q|^{2d-2}
\left(\frac{\Ima_I(q^d)}{y}\right)^2>0.
\end{equation}
\end{proof}

The power map is proper because $|q^d|=|q|^d$. Its Brouwer degree is $d$ by Lemma~\ref{lem:power} and the regular-value formula. Real target powers can be obtained by choosing a polar argument and dividing it by $d$, as described in Subsection~\ref{subsec:real-targets}, so the power map is surjective on all of $\HH$.

\section{The scalar radius and the local degree}\label{sec:reduction}
We keep the notation of Subsection~\ref{subsec:principal}, with a non-real target $b=a+I\rho$ and the polynomial $R_b$ of~\eqref{eq:radial}.

\begin{theorem}[Scalar reduction and a sharp bound]\label{thm:count}
Let $d=2m\geq2$ and $\deg\phi=r\geq1$. For every non-real $b$, the equation $P_{d,\phi}(q)=b$ has precisely
\[
d\,\#\{t>0\mid R_b(t)=0\}
\]
distinct quaternionic solutions, where the positive roots of $R_b$ are counted without multiplicity. All lie in $\CC_I$. For each such $t$, they are the $d$ complex roots in $\CC_I$ of
\begin{equation}\label{eq:lift}
q^d=a-\phi(t)+I\rho.
\end{equation}
Every positive root $t$ of $R_b$ satisfies $t\geq\rho$. The value $t=\rho$ is a root exactly when $\phi(\rho)=a$. The polynomial $R_b$ is non-zero of degree at most $2r$, giving at most $2dr$ solutions. In the affine case $\phi(t)=\pm t+a$, one has $R_b\equiv\rho^2>0$ and no solutions. A quaternionic solution over $t$, that is, with $|q|^d=t$, is simple exactly when $R_b'(t)\neq0$. Its local degree is
\begin{equation}\label{eq:sign}
-\sgn R_b'(t).
\end{equation}
If $\phi(\rho)=a$, then $R_b'(\rho)=-2\rho<0$, so this radial root is simple and its $d$ quaternionic preimages all have local degree $+1$. If $R_b'(t)=0$ and $R_b''(t)\neq0$ at a positive root of $R_b$, then all $d$ preimages over that radius are ordinary fold critical points, as explained in Remark~\ref{rem:radial-contact}.

For every pair $(d,r)$ in the stated range, there are $\phi$ and a non-real target for which all $2dr$ roots occur and are simple.
\end{theorem}

\begin{proof}[Proof of Theorem~\ref{thm:count}, apart from sharpness]
Suppose $P_{d,\phi}(q)=b$ and set $t=|q|^d$. Since $b$ is non-real, $q\neq0$ and
\[
q^d=b-\phi(t).
\]
The right side is non-real. Since $q$ commutes with its own power and $\phi(t)$ is real, it commutes with $b$. Writing $q=x+v$, the identity $qI-Iq=2v\times I=0$ gives $v\in\RR I$, hence $q\in\CC_I$. Taking norms gives $R_b(t)=0$.

Conversely, let $t>0$ satisfy $R_b(t)=0$. Then $|b-\phi(t)|=t$, and Lemma~\ref{lem:power} gives exactly $d$ roots of~\eqref{eq:lift}. Each has $|q|^d=t$ and therefore solves the original equation. Different radii give disjoint sets of roots. The equality $t^2=(a-\phi(t))^2+\rho^2$ also gives $t\geq\rho$.

To determine the derivative, factor the map as
\begin{equation}\label{eq:factor}
P_{d,\phi}=T_\phi\circ Q_d,\qquad T_\phi(w)=w+\phi(|w|).
\end{equation}
At $w\neq0$, with $t=|w|$, the radial term contributes a rank-one derivative. The determinant is
\[
\det DT_\phi(w)=1+\phi'(t)\frac{\Rea w}{t}.
\]
At a preimage of $b$, this becomes
\begin{equation}\label{eq:radialjac}
\det DT_\phi(w)
=1+\frac{(a-\phi(t))\phi'(t)}{t}
=-\frac{R_b'(t)}{2t}.
\end{equation}
Combining this with~\eqref{eq:powerjac} proves simplicity and the local degree formula.

For $r\geq2$, the leading coefficient of $R_b$ is the square of the leading coefficient of $\phi$, so its degree is $2r$. For $r=1$, the degree is at most two. It cannot vanish identically. Indeed, such vanishing would require the slope of $\phi$ to be $+1$ or $-1$, its constant term to be $a$, and $\rho^2=0$, contrary to the non-real target assumption. If $R_b$ is constant, then it is non-zero and has no roots. The root bound follows in all cases.
\end{proof}

The real coordinate degree of $P_{d,\phi}$ is exactly $N=dr$. For $r>1$, the leading radial term has degree $dr>d$. For $r=1$, the imaginary part of $q^d$ is a non-zero homogeneous polynomial map of degree $d$, so the real radial term cannot cancel it. Thus the sharp bound is $2N$ isolated roots for a degree-$N$ map in this class.

\begin{remark}[Degrees of the proper radial maps]\label{rem:radial-degree}
In the affine case write $\phi(t)=\lambda t+c$. The map $P_{d,\phi}$ is proper for $\deg\phi\geq2$, and in the affine case exactly when $|\lambda|\neq1$. For these maps
\begin{equation}\label{eq:radial-brouwer-degree}
\deg P_{d,\phi}=\begin{cases}
d,& \phi(t)=\lambda t+c,\quad |\lambda|<1,\\
0,& \phi(t)=\lambda t+c,\quad |\lambda|>1,\\
0,& \deg\phi\geq2.
\end{cases}
\end{equation}
Properness follows from $|P_{d,\phi}(q)|\geq\bigl||\phi(t)|-t\bigr|$ with $t=|q|^d$. For $\lambda=\pm1$, the rays with $q^d=-\lambda|q|^d$ map to $c$. At a non-real regular target, $R_b(0)>0$, so the local signs~\eqref{eq:sign} alternate, beginning with $+1$. Their sum is one when $|\lambda|<1$, since $R_b(t)\to-\infty$, and zero in the other proper cases, since $R_b(t)\to+\infty$. Each radius has $d$ preimages.
\end{remark}

The scalar equation also describes degeneracy. A multiple positive root of $R_b$ produces non-simple quaternionic preimages. Thus root counts and the signs of their Jacobians can be read from a single real polynomial without assigning a complex multiplicity to a general two-sided expression.

\begin{remark}[Radial contact and branching]\label{rem:radial-contact}
Put $w_b(t)=b-\phi(t)$. The equation $R_b(t)=0$ means that the moving point $w_b(t)$ lies on the sphere of radius $t$ centered at zero. Equivalently, the curve $t\mapsto(t,w_b(t))$ meets the smooth cone $|w|^2=t^2$ in $(0,\infty)\times\HH$. At a root, $R_b'(t)=0$ means that this curve is tangent to the cone. If also $R_b''(t)\neq0$, then the contact is quadratic and the corresponding quaternionic preimages are ordinary folds. Indeed, writing $w=u+v$ gives $T_\phi(u,v)=(u+\phi(\sqrt{u^2+|v|^2}),v)$. At a critical preimage, put $u=a-\phi(t)$. The identity $R_b'(t)=-2u\phi'(t)-2t=0$ gives $u\phi'(t)=-t<0$, since $t\geq\rho>0$. Thus $u\neq0$, and the vanishing first derivative gives
\[
1+\frac{u}{t}\phi'(t)=0
\qquad\Longrightarrow\qquad
\phi'(t)=-\frac{t}{u}.
\]
For the partial derivatives with respect to $u$, the imaginary coordinate $v$ is held fixed, so $|v|=\rho$ remains constant. Substitution into the second derivative at the critical preimage yields
\[
\frac{\partial^2}{\partial u^2}\bigl(u+\phi(\sqrt{u^2+|v|^2})\bigr)
=-\frac{u}{2t^2}R_b''(t).
\]
Since $T_\phi$ has the form $(u,v)\mapsto(f(u,v),v)$, its Jacobian determinant is $\partial f/\partial u$, and its critical set is defined locally by $\partial f/\partial u=0$. At a critical point, the kernel of $DT_\phi$ is spanned by $\partial/\partial u$. The derivative of $\det DT_\phi$ along this kernel is the second partial derivative just computed. When it is non-zero, the critical set is a smooth hypersurface transverse to the kernel, which is the corank-one fold criterion \cite[Chapter~III, \S4]{golubitskyguillemin}. A multiple radial root alone permits higher degeneracy, so tangency by itself does not imply a fold.
\end{remark}

\subsection{Real targets and spherical fibers}\label{subsec:real-targets}
Let $b=a\in\RR$. For $t>0$, the radial condition reduces to
\[
\phi(t)+t=a
\quad\text{or}\quad
\phi(t)-t=a.
\]
In the first case $q^d=t$, and in the second case $q^d=-t$. For $d=2m$, the first fiber consists of two real points and $m-1$ conjugacy two-spheres, all at radius $t^{1/d}$. The second consists of $m$ conjugacy two-spheres. This follows from the polar expression $q=t^{1/d}(\cos\theta+J\sin\theta)$ with $J^2=-1$ and the conditions $d\theta\in2\pi\ZZ$ or $d\theta\in\pi+2\pi\ZZ$. The zero quaternion is a solution exactly when $\phi(0)=a$.

The conjugacy spheres correspond to spherical roots in the terminology of the one-sided and slice regular theories \cite{pogoruishapiro,gentilistoppatostruppa}, while their radial locations here are determined by $\phi(t)\pm t=a$. These fibers explain the non-real target assumption in the finite root bound. A real radial equation can also vanish identically, as at the exceptional parameter in Theorem~\ref{thm:transition}.

\section{Sharpness and an explicit quartic}\label{sec:sharp}
\begin{proof}[Proof of sharpness in Theorem~\ref{thm:count}]
Fix $\rho>0$ and choose distinct numbers $\rho<t_1<\cdots<t_r$. Let
\[
g(t)=\prod_{j=1}^r(t-t_j),\qquad \phi_L(t)=L g(t),\qquad b=I\rho.
\]
The scalar polynomial is
\[
R_b(t)=L^2g(t)^2+\rho^2-t^2.
\]
At each $t_j$ its value is negative. Choose numbers $u_0,\ldots,u_r$ interlacing the $t_j$ with $\rho<u_0<t_1$ and $t_r<u_r$. At all these points $g(u_j)\neq0$. Choose $L>0$ so that
\[
L^2>\max_{0\leq j\leq r}\frac{u_j^2-\rho^2}{g(u_j)^2}.
\]
Then every $R_b(u_j)$ is positive. There is consequently a root in each of the $2r$ disjoint intervals between consecutive points of
\[
u_0<t_1<u_1<t_2<\cdots<u_{r-1}<t_r<u_r.
\]
Take $L>1$ when $r=1$, so the quadratic coefficient is $L^2-1>0$. The polynomial has degree exactly $2r$, so these are all its roots. They are distinct, positive, and simple. Theorem~\ref{thm:count} gives exactly $2dr$ simple quaternionic roots. The signs of $R_b'$ alternate, yielding $dr$ roots of each local degree.
\end{proof}

\begin{proposition}[A quadratic power with a quartic radial term]\label{prop:quartic}
The two-sided quaternionic equation
\begin{equation}\label{eq:quartic}
q^2+4(|q|^2-2)(|q|^2-4)-i=0
\end{equation}
has exactly eight roots, all simple and all in $\CC_i$. Four have local degree $+1$ and four have local degree $-1$.
\end{proposition}
\begin{proof}
Here $d=r=2$ and the radial polynomial is
\[
R(t)=16(t-2)^2(t-4)^2+1-t^2.
\]
Its values at $1,2,3,4,5$ are respectively $144,-3,8,-15,120$. Hence its four roots lie one each in $(1,2)$, $(2,3)$, $(3,4)$, and $(4,5)$. A quartic with four distinct real roots has no further roots and every root is simple. Each radius gives two complex square roots of $i-4(t-2)(t-4)$. Formula~\eqref{eq:sign} gives the asserted signs.
\end{proof}

The radial polynomial factors as
\begin{equation}\label{eq:quartic-factorization}
R(t)=(4t^2-21t+25)(4t^2-27t+41),
\end{equation}
so its four positive roots are $(21\pm\sqrt{41})/8$ and $(27\pm\sqrt{73})/8$. Each value gives the two square roots in $\CC_i$ of $i-4(t-2)(t-4)$, while the sign argument above proves simplicity. The polynomial degree is four in the real coordinates. For a unilateral polynomial of degree $N$, the number of distinct isolated roots plus twice the number of distinct spherical roots is at most $N$ \cite{pogoruishapiro}. The eight isolated roots in~\eqref{eq:quartic} exceed that unilateral bound of four.

\section{The linear radial term and the loss of degree}\label{sec:transition}
\begin{theorem}[A degree transition and its image]\label{thm:transition}
Let $d=2m\geq2$ and $\lambda\geq0$. Put
\[
P_\lambda(q)=q^d+\lambda|q|^d.
\]
For $0\leq\lambda<1$, the map is proper, surjective, and has Brouwer degree $d$. Every non-real target has exactly $d$ simple preimages, all of local degree $+1$.

For $\lambda>1$, the map is proper of Brouwer degree zero and has image
\begin{equation}\label{eq:cone}
K_\lambda=\{a+v\in\HH\mid a\geq\sqrt{\lambda^2-1}\,|v|\}.
\end{equation}
Every non-real target in the interior of this cone has $2d$ simple preimages, with $d$ of each local degree. Every non-zero boundary target is non-real and has exactly $d$ preimages, all ordinary fold singularities. Targets outside the cone have no preimages.

At $\lambda=1$, the image is $\{b\mid\Rea b>0\}\cup\{0\}$, and the map is not proper. Each non-real target of positive real part has exactly $d$ simple preimages.
\end{theorem}

For $\phi(t)=\lambda t$, equation~\eqref{eq:radial} becomes
\begin{equation}\label{eq:quadratic}
(\lambda^2-1)t^2-2\lambda a t+a^2+\rho^2=0.
\end{equation}
This quadratic displays the transition directly. For $\lambda<0$, the image is the reflection of that for $|\lambda|$ through the origin, replacing $a$ by $-a$ in~\eqref{eq:cone} when $\lambda<-1$ and using $u^d=1$ for the zero-fiber construction at $\lambda=-1$, consistently with Remark~\ref{rem:radial-degree}.

\begin{proof}[Proof of Theorem~\ref{thm:transition}]
For $0\leq\lambda<1$, the constant term of~\eqref{eq:quadratic} is positive and its leading coefficient is negative. Its two roots have opposite signs, so there is exactly one positive root. At this root $R_b'<0$. Theorem~\ref{thm:count} gives $d$ simple preimages of local degree $+1$. Moreover
\[
|P_\lambda(q)|\geq(1-\lambda)|q|^d.
\]
The map is proper, its Brouwer degree is $d$, and non-zero degree gives surjectivity.

For $\lambda>1$, put $\Delta=a^2-(\lambda^2-1)\rho^2$, the reduced discriminant of~\eqref{eq:quadratic}. The standard quadratic discriminant is $4\Delta$. There are two distinct positive roots exactly when $a>\sqrt{\lambda^2-1}\rho$, and then
\begin{equation}\label{eq:radii}
t_\pm=\frac{\lambda a\pm\sqrt\Delta}{\lambda^2-1}.
\end{equation}
Their derivatives satisfy $R_b'(t_+)=2\sqrt\Delta$ and $R_b'(t_-)=-2\sqrt\Delta$. By~\eqref{eq:sign}, the smaller radius $t_-$ yields $d$ preimages of local degree $+1$, while the larger radius $t_+$ yields $d$ preimages of local degree $-1$. At equality there is one positive double root, producing $d$ preimages. Below this boundary there are no positive roots. Real targets follow directly from $q^d=\pm|q|^d$. Every positive real target is attained, zero is attained only at zero, and negative real targets are omitted. This proves the full cone image~\eqref{eq:cone}.

The estimate $|P_\lambda(q)|\geq(\lambda-1)|q|^d$ proves properness. A target outside the image has an empty preimage, so the Brouwer degree is zero. The fold assertion is proved just below.

For $\lambda=1$, the scalar equation is $-2at+a^2+\rho^2=0$. For a non-real target it has a positive solution exactly when $a>0$, and that root is simple. The real part of $w+|w|$ is non-negative and vanishes only when $w$ is non-positive real, in which case $w+|w|=0$. This proves the asserted image, using surjectivity of $Q_d$. Choose a unit quaternion $u$ with $u^d=-1$. Then $P_1(ru)=0$ for every $r\geq0$, which proves failure of properness. The family $(\lambda,q)\mapsto P_\lambda(q)$ is not proper on any parameter interval containing $1$, since the preimage of zero contains this unbounded fiber. Proper homotopy invariance therefore does not apply across the threshold.
\end{proof}

\subsection{The apex and the threshold zero fiber}\label{subsec:apex-threshold}
For $\lambda>1$, the cone apex has the unique preimage $q=0$, since $|P_\lambda(q)|\geq(\lambda-1)|q|^d$. The map is homogeneous of degree $d\geq2$, so $DP_\lambda(0)=0$. Its derivative has corank four at the origin, which excludes the ordinary fold normal form. Thus the apex has a homogeneous degeneracy distinct from the corank-one folds over the non-zero boundary.

\begin{figure}[htbp]
\centering
\begin{tikzpicture}[>={Stealth[length=2mm]},font=\small]
\begin{scope}
\draw[->] (-2.4,0) -- (2.6,0) node[below] {$\Rea q$};
\draw[->] (0,-2.4) -- (0,2.6) node[left] {$\Ima_J q$};
\draw[dashed] (0,0) circle (1);
\foreach \a in {30,90,150,210,270,330} {
  \draw[very thick] (0,0) -- (\a:2.2);
  \fill (\a:1) circle (1.6pt);
}
\draw (0.7,0) arc (0:30:0.7);
\node at (15:1.3) {$\theta_0$};
\node at (30:2.55) {$\Sigma_0$};
\node at (330:2.55) {$\Sigma_0$};
\node at (150:2.55) {$\Sigma_2$};
\node at (210:2.55) {$\Sigma_2$};
\node at (0.4,2.15) {$\Sigma_1$};
\node at (0.4,-2.15) {$\Sigma_1$};
\fill (0,0) circle (1.6pt);
\node at (0,-3.0) {(a)};
\end{scope}
\begin{scope}[xshift=6.6cm]
\draw[->] (-2.9,0) -- (3.0,0) node[below] {$\RR$};
\fill[gray!15] (0,0) -- (2.08,1.2) arc (90:-90:0.3 and 1.2) -- cycle;
\fill[gray!15] (0,0) -- (-2.08,1.2) arc (90:270:0.3 and 1.2) -- cycle;
\draw (0,0) -- (2.08,1.2);
\draw (0,0) -- (2.08,-1.2);
\draw (2.08,0) ellipse (0.3 and 1.2);
\draw (0,0) -- (-2.08,1.2);
\draw (0,0) -- (-2.08,-1.2);
\draw (-2.08,0) ellipse (0.3 and 1.2);
\draw[fill=gray!8,fill opacity=0.6] (0,0) ellipse (0.6 and 2.4);
\fill (0,0) circle (1.6pt);
\node at (2.08,1.5) {$\Sigma_0$};
\node at (-2.08,1.5) {$\Sigma_2$};
\node at (0,2.65) {$\Sigma_1$};
\node at (0,-3.0) {(b)};
\end{scope}
\end{tikzpicture}
\caption{Threshold zero fiber~\eqref{eq:threshold-zero-fiber} for $d=6$. In (a), its intersection with $\CC_J$ consists of six rays through the roots of $u^6=-1$ on the unit circle. Each label $\Sigma_j$ in (a) identifies a ray through a point of $\Sigma_j\cap\CC_J$. In (b), one imaginary dimension is suppressed. The two-sphere links $\Sigma_j$ appear as circles, and their cones occupy angles $\pi/6$, $\pi/2$ and $5\pi/6$ from the positive real axis. The middle component away from the origin is $\Ima\HH\setminus\{0\}$. This imaginary component occurs exactly when $m=d/2$ is odd.}
\label{fig:threshold-zero-fiber}
\end{figure}
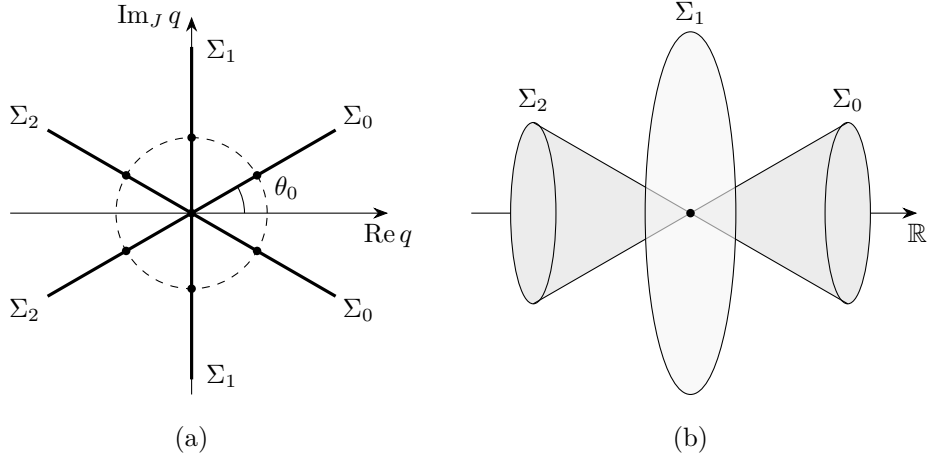

At $\lambda=1$, every positive real target $a$ has the two real preimages $\pm(a/2)^{1/d}$ and $m-1$ conjugacy two-spheres, all of radius $(a/2)^{1/d}$, as follows from $q^d=|q|^d=a/2$ and Subsection~\ref{subsec:real-targets}. The zero fiber can be described explicitly. Set
\[
\theta_j=\frac{(2j+1)\pi}{d},\qquad 0\leq j\leq m-1,\qquad d=2m
\]
and define the disjoint conjugacy two-spheres
\[
\Sigma_j=\{\cos\theta_j+J\sin\theta_j\mid J\in\Ima\HH,\ |J|=1\}
\subset\mathbb S^3.
\]
Every non-zero quaternion has a unique polar angle in $[0,\pi]$. The condition $u^d=-1$ for a unit quaternion excludes the endpoints and gives exactly the angles $\theta_j$. Consequently
\begin{equation}\label{eq:threshold-zero-fiber}
P_1^{-1}(0)=\{0\}\cup\bigcup_{j=0}^{m-1}\{ru\mid r>0,\ u\in\Sigma_j\}.
\end{equation}
This is the cone over $m$ disjoint two-spheres, illustrated for $d=6$ in Figure~\ref{fig:threshold-zero-fiber}. Each component away from the origin is diffeomorphic to $(0,\infty)\times\mathbb S^2$, so the zero fiber has real dimension three. For $d=2$, it is the imaginary hyperplane $\Ima\HH$. Formula~\eqref{eq:threshold-zero-fiber} describes the unbounded fiber responsible for the failure of properness at the degree transition.

\subsection{Non-zero boundary points and ordinary folds}\label{subsec:folds}
Fix $\lambda>1$. A boundary point with $v=0$ has $a=0$, so every non-zero boundary point is non-real. Factor $P_\lambda=T_\lambda\circ Q_d$, where $T_\lambda(w)=w+\lambda|w|$. Write $w=u+v$, with $u\in\RR$ and $v\in\RR^3$. Then
\[
T_\lambda(u,v)=\bigl(u+\lambda\sqrt{u^2+|v|^2}\,,\ v\bigr).
\]
At a preimage of a non-zero cone boundary point, $v\neq0$. The derivative of the first coordinate with respect to $u$ vanishes exactly at
\[
u=-\frac{|v|}{\sqrt{\lambda^2-1}}.
\]
At this critical point, Remark~\ref{rem:radial-contact} applies with $\phi(t)=\lambda t$ and $t=\sqrt{u^2+|v|^2}$. At this point $u=-t/\lambda$ and $t^2=\lambda^2|v|^2/(\lambda^2-1)$. Since $R_b''(t)=2(\lambda^2-1)>0$, the second derivative is
\[
-\frac{u}{2t^2}R_b''(t)
=\frac{\lambda^2-1}{\lambda t}
=\frac{(\lambda^2-1)^{3/2}}{\lambda^2|v|}
=\frac{\lambda|v|^2}{(u^2+|v|^2)^{3/2}}>0.
\]
The last three coordinates are unchanged. The critical manifold $S$ is defined locally by $h(u,v)=u+|v|/\sqrt{\lambda^2-1}=0$. Along $S$, the kernel of $DT_\lambda$ is spanned by $\partial/\partial u$, and $dh(\partial/\partial u)=1$. Thus this kernel is transverse to $TS$. The fold normal form \cite[Chapter~III, \S4]{golubitskyguillemin} therefore gives the local form $(x_1,x_2,x_3,x_4)\mapsto(x_1^2,x_2,x_3,x_4)$ after smooth coordinate changes. By Lemma~\ref{lem:power}, $Q_d$ is an orientation-preserving local diffeomorphism at each preimage of the non-real value $w$. The same fold description and local-degree signs therefore hold at all $d$ quaternionic lifts.

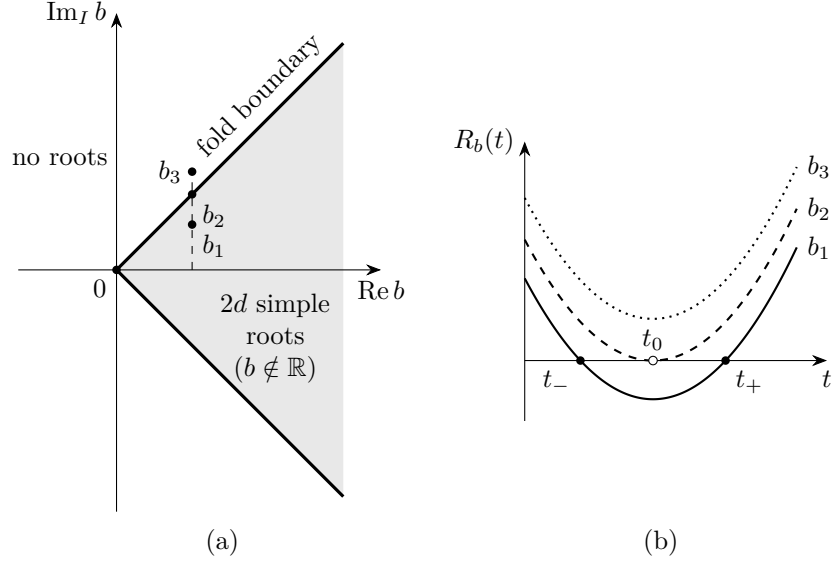
\begin{figure}[htbp]
\centering
\begin{tikzpicture}[>={Stealth[length=2mm]},font=\small]
\begin{scope}
\fill[gray!18] (0,0) -- (3,3) -- (3,-3) -- cycle;
\draw[->] (-1.3,0) -- (3.5,0) node[below] {$\Rea b$};
\draw[->] (0,-3.2) -- (0,3.4) node[left] {$\Ima_I b$};
\draw[very thick] (0,0) -- (3,3);
\draw[very thick] (0,0) -- (3,-3);
\node[rotate=45,above] at (2.1,2.1) {fold boundary};
\node[align=center] at (2.1,-0.9) {$2d$ simple\\roots\\($b\notin\RR$)};
\node at (-0.75,1.5) {no roots};
\draw[dashed] (1,0) -- (1,1.3);
\fill (1,0.6) circle (1.6pt) node[below right] {$b_1$};
\fill (1,1) circle (1.6pt) node[below right] {$b_2$};
\fill (1,1.3) circle (1.6pt) node[left] {$b_3$};
\fill (0,0) circle (1.6pt) node[below left] {$0$};
\node at (1.4,-3.6) {(a)};
\end{scope}
\begin{scope}[xshift=5.4cm,yshift=-1.2cm]
\draw[->] (0,0) -- (4.0,0) node[below] {$t$};
\draw[->] (0,-0.8) -- (0,2.9) node[left] {$R_b(t)$};
\draw[thick,domain=0:3,samples=80,smooth,variable=\t] plot ({1.2*\t},{0.8*(\t*\t-2*sqrt(2)*\t+1.36)}) node[right] {$b_1$};
\draw[thick,dashed,domain=0:3,samples=80,smooth,variable=\t] plot ({1.2*\t},{0.8*(\t*\t-2*sqrt(2)*\t+2)}) node[right] {$b_2$};
\draw[thick,dotted,domain=0:3,samples=80,smooth,variable=\t] plot ({1.2*\t},{0.8*(\t*\t-2*sqrt(2)*\t+2.69)}) node[right] {$b_3$};
\fill ({1.2*(sqrt(2)-0.8)},0) circle (1.6pt) node[below left] {$t_-$};
\fill ({1.2*(sqrt(2)+0.8)},0) circle (1.6pt) node[below right] {$t_+$};
\draw[fill=white] ({1.2*sqrt(2)},0) circle (1.6pt) node[above] {$t_0$};
\node at (1.8,-2.4) {(b)};
\end{scope}
\end{tikzpicture}
\caption{Cone section and radial fibers for $\lambda=\sqrt2$. In (a), $b_1=1+0.6I$, $b_2=1+I$ and $b_3=1+1.3I$ lie inside, on and outside the cone~\eqref{eq:cone}. Interior non-real targets have $2d$ simple preimages, and non-zero boundary targets have $d$ ordinary folds. Positive real targets have spherical fibers. In (b), the radial polynomials~\eqref{eq:quadratic} have two positive roots for $b_1$, a double root $t_0=\sqrt2$ for $b_2$ and no real root for $b_3$. Each smaller-radius lift has local degree $+1$ and each larger-radius lift $-1$.}
\label{fig:cone-folds}
\end{figure}

As a target crosses the cone boundary, each fold contributes a pair of simple roots with opposite local degree on the interior side and no roots on the exterior side. The total signed count remains zero. Figure~\ref{fig:cone-folds} compares a cone section with the radial polynomials for three targets. The two positive radii~\eqref{eq:radii} merge at the boundary into the double root that produces the folds.

\subsection{Lower-degree perturbations on the non-zero-degree side}\label{subsec:lower-degree}
\begin{corollary}\label{cor:lower}
Let $0\leq\lambda<1$, and let $E:\HH\to\HH$ be any two-sided polynomial map of degree less than $d$. Then $P_\lambda+E$ is proper of Brouwer degree $d$, hence surjective. Every regular target has at least $d$ preimages.
\end{corollary}
\begin{proof}
The homotopy $P_\lambda+sE$, for $0\leq s\leq1$, satisfies
\[
|P_\lambda(q)+sE(q)|\geq(1-\lambda)|q|^d-O(|q|^{d-1}+1)
\]
uniformly in $s$. Thus it is a proper homotopy and preserves degree. The regular-value formula gives the last assertion.
\end{proof}

Properness fails at the intervening parameter $\lambda=1$, so the change of degree does not violate homotopy invariance. On the degree-zero side, root existence depends on the target's position relative to the cone.

\section{Central perturbations and planar fibers}\label{sec:central}
A real-valued perturbation preserves the commutation argument even when it depends on more than the norm. The planar description below isolates the contribution of centrality.

\begin{proposition}[Planar localization]\label{prop:central}
Let $d\geq1$, let $f:\HH\to\RR$ be a real polynomial, and let $P(q)=q^d+f(q)$ have real coordinate degree $N$. For a non-real target $b=a+I\rho$ with $\rho>0$, every solution lies in $\CC_I$. Define real polynomials $A_I,B_d$ by
\[
A_I(x,y)=\Rea\bigl((x+Iy)^d\bigr)+f(x+Iy),\qquad B_d(x,y)=\Ima_I\bigl((x+Iy)^d\bigr).
\]
The fiber is identified with the real solutions of
\begin{equation}\label{eq:central-planar}
A_I(x,y)=a,\qquad B_d(x,y)=\rho.
\end{equation}
At such a solution, $y\neq0$ and
\begin{equation}\label{eq:central-jacobian}
\det DP(x+Iy)=\left(\frac{\rho}{y}\right)^2
\det D(A_I,B_d)(x,y).
\end{equation}
There are at most $dN$ simple quaternionic solutions. If the complexification of~\eqref{eq:central-planar} has a finite zero set, then the same bound holds for all distinct solutions.
\end{proposition}
\begin{proof}
The identity $q^d=b-f(q)$ implies $qb=bq$, giving $q\in\CC_I$. A real $q$ cannot map to a non-real target, so $y\neq0$. The coordinate equations follow by restriction to this slice.

Use the same basis $(1,I,J,IJ)$ in the domain and codomain, with $J$ a unit imaginary quaternion orthogonal to $I$. The calculation in Lemma~\ref{lem:power} gives transverse scalar multiplication by $\rho/y$ for $DQ_d$. The derivative of the real-valued term takes values in the first coordinate. With respect to the splitting $\HH=\CC_I\oplus\CC_I^\perp$, the slice $\CC_I$ is therefore preserved by $DP$, and the lower-left $2\times2$ block vanishes. On $\CC_I^\perp$, the derivative $DQ_d$ acts as multiplication by $\rho/y$, while $Df$ takes real values and contributes to the upper-right block alone. Thus $DP$ is block upper triangular, with diagonal blocks $D(A_I,B_d)$ and $(\rho/y)\operatorname{Id}_{\RR^2}$. This proves~\eqref{eq:central-jacobian}.

The planar equations have degrees at most $N$ and $d$. A simple quaternionic solution has an invertible planar Jacobian and is therefore an isolated zero of the complexified equations. Isolated zeros lie off any common component, so removing the greatest common divisor leaves them in a coprime system to which the planar B\'ezout bound $dN$ applies \cite{fulton}. Under the finite-zero-set hypothesis it also bounds all distinct real solutions. For each fixed non-real target and each fixed degree bound on $f$, the finite-zero-set hypothesis holds on an open dense subset of the full real coefficient space. Indeed, failure requires $A_I-a$ to vanish identically on an irreducible component of $B_d-\rho=0$. On each of the finitely many components, this condition defines either an empty set or a proper affine subspace of coefficients, since it is affine in $f$ and changing the constant term destroys the vanishing.
\end{proof}

\begin{proposition}[Twelve roots for a central quartic]\label{prop:central-quartic}
Put $p(x)=(x-1)(x-2)(x-3)(x-4)$. The map
\begin{equation}\label{eq:central-quartic}
P(q)=q^4-\Rea(q^4)+p(\Rea q)
\end{equation}
is a real-valued perturbation of $q^4$ of real coordinate degree four. The equation $P(q)=i$ has exactly twelve roots, all simple, with six of each local degree.
\end{proposition}
\begin{proof}
Proposition~\ref{prop:central} places every root in $\CC_i$. Writing $q=x+iy$ gives
\[
p(x)=0,\qquad 4xy(x^2-y^2)=1.
\]
For each $x=j\in\{1,2,3,4\}$, the second equation is
\[
y^3-j^2y+\frac{1}{4j}=0.
\]
Its discriminant is $4j^6-27/(16j^2)>0$, since $64j^8\geq64>27$. It has exactly three distinct real roots, all non-zero. These give twelve distinct quaternionic roots. The planar Jacobian determinant at a root is $4j\,p'(j)(j^2-3y^2)$, which is non-zero because both univariate roots are simple. Formula~\eqref{eq:central-jacobian} proves quaternionic simplicity.

For fixed $j$, write the three roots in increasing order as $y_1,y_2,y_3$. The critical points of the cubic separate them as
\[
y_1<-\frac{j}{\sqrt3}<y_2<\frac{j}{\sqrt3}<y_3.
\]
The factors $j^2-3y_\ell^2$ therefore have signs $-1,+1,-1$. The local degrees over $x=j$ sum to $-\sgn p'(j)$. Since $\sgn p'(j)=(-1)^{4-j}$, the total signed count is
\[
-\sum_{j=1}^4(-1)^{4-j}=0.
\]
All twelve roots are simple, so exactly six have each local degree.
\end{proof}

This example lies even in the conjugation-invariant central class $f(q)=\psi(\Rea q,|q|^2)$. Indeed, writing $x=\Rea q$ and $s=|q|^2$ gives $\Rea(q^4)=8x^4-8x^2s+s^2$. Its twelve roots exceed the radial bound $2N=8$ at the same coordinate degree. Slice localization and conjugation symmetry therefore persist beyond the family $P_{d,\phi}$. Dependence on $x$ retains angular information on a slice. Proposition~\ref{prop:norm-square} gives a scalar equation in $s$ when the perturbation depends on $s$ alone, and it reaches sixteen roots at this degree.

The map~\eqref{eq:central-quartic} is not coercive. Every purely imaginary $q$ satisfies $q^4=|q|^4\in\RR$, so $P$ takes the constant value $p(0)=24$ on the whole of $\Ima\HH$. Its twelve simple roots persist under small perturbations by the implicit function theorem. Lemma~\ref{lem:stability} does not apply to it, however, and further roots may enter from infinity under perturbation. The radial examples of Proposition~\ref{prop:norm-square} are coercive, and Section~\ref{sec:stability} shows that their counts persist on an open set.

Within the family $q^4-\Rea(q^4)+p(\Rea q)$ with a non-zero polynomial $p$ of degree at most four, twelve is the largest number of distinct roots over $i$. Indeed, the equation $p(x)=0$ allows at most four real values of $x$, the value $x=0$ gives no solution, and each non-zero value gives a cubic equation in $y$ with at most three roots.

\subsection{A radius measured by the squared norm}\label{subsec:norm-square}
The bound $2N$ of Theorem~\ref{thm:count} comes from the dependence of the radial term on $|q|^d$. A radial term in $|q|^2$ keeps the factorization through the power map and the one-radius equation. Proposition~\ref{prop:norm-square} attains the planar bound $dN$ of Proposition~\ref{prop:central} for every even $N\geq d$ within this conjugation-invariant radial class. The construction also allows odd exponents $d$.

\begin{proposition}[A radius measured by the squared norm]\label{prop:norm-square}
Let $d\geq1$, let $\psi\in\RR[u]$ have degree $k\geq1$, and put $P(q)=q^d+\psi(|q|^2)$. Its real coordinate degree is $N=\max\{d,2k\}$. For a non-real target $b=a+I\rho$ with $\rho>0$, put
\begin{equation}\label{eq:norm-square-radial}
S_b(u)=(a-\psi(u))^2+\rho^2-u^d.
\end{equation}
The equation $P(q)=b$ has exactly $d\,\#\{u>0\mid S_b(u)=0\}$ distinct solutions, all in $\CC_I$. For each such $u$, they are the $d$ roots in $\CC_I$ of $q^d=a-\psi(u)+I\rho$. The polynomial $S_b$ is non-zero of degree at most $N$, so there are at most $dN$ solutions. If $2k<d$, then there are at most $d(2k+1)$ solutions. A solution over $u$ is simple exactly when $S_b'(u)\neq0$, and its local degree is then $-\sgn S_b'(u)$. If $2k\geq d$, then there are $\psi$ and $b$ for which all $dN$ solutions occur and are simple, half of each local degree.
\end{proposition}
\begin{proof}
For the degree, the leading term of $\psi(|q|^2)$ is a non-zero multiple of $|q|^{2k}$. When $2k=d\geq2$, the imaginary part of $q^d$ is a non-zero homogeneous map of degree $d$ that this real term cannot cancel. Hence $N=\max\{d,2k\}$ in all cases.

Suppose that $P(q)=b$ and put $u=|q|^2$. The quaternion $q^d=b-\psi(u)$ is non-real, so $q\neq0$ and $q$ commutes with $b$, which places $q$ in $\CC_I$. Taking norms gives $u^d=|q|^{2d}=(a-\psi(u))^2+\rho^2$, that is, $S_b(u)=0$. Conversely, let $u>0$ satisfy $S_b(u)=0$ and put $w=a-\psi(u)+I\rho$. Then $|w|=u^{d/2}$, and Lemma~\ref{lem:power} gives $d$ roots of $q^d=w$ in $\CC_I$, each with $|q|^2=u$ and hence $P(q)=w+\psi(u)=b$. Different values of $u$ give disjoint sets of roots.

For the derivative, write $P=T\circ Q_d$ with $T(w)=w+\varphi(|w|)$ and $\varphi(t)=\psi(t^{2/d})$, smooth for $t>0$. The computation of~\eqref{eq:radialjac} uses only the chain rule, so it applies to $\varphi$ and gives $\det DT(w)=-R'(t)/(2t)$ with $R(t)=(a-\varphi(t))^2+\rho^2-t^2=S_b(t^{2/d})$. The substitution $u=t^{2/d}$ has positive derivative, so $\sgn R'(t)=\sgn S_b'(u)$. Together with~\eqref{eq:powerjac}, this proves the assertions on simplicity and local degree.

The value $S_b(0)=(a-\psi(0))^2+\rho^2$ is positive, so $S_b$ is non-zero. With $c$ the leading coefficient of $\psi$, the leading coefficient of $S_b$ is $c^2$ when $2k>d$, it is $-1$ when $2k<d$, and its degree is at most $d$ when $2k=d$. The degree of $S_b$ is thus at most $N$, which gives the bound $dN$. When $2k<d$, the only term above degree $2k$ is $-u^d$. There are at most $2k+1$ sign changes in the coefficient sequence, so Descartes' rule gives at most $2k+1$ positive roots, counted with multiplicity, and hence at most $d(2k+1)$ quaternionic solutions. For even $d$, this bound is strictly below $d^2$.

For sharpness, let $2k\geq d$, so that $N=2k$. Fix $\rho>0$, put $b=I\rho$, and choose numbers $\rho^{2/d}<u_1<\cdots<u_k$ and $v_0,\ldots,v_k$ with
\[
\rho^{2/d}<v_0<u_1<v_1<\cdots<u_k<v_k.
\]
Put $g(u)=\prod_{j=1}^k(u-u_j)$ and $\psi=Lg$ with $L>0$, so that $S_b(u)=L^2g(u)^2+\rho^2-u^d$. At each $u_j$ its value $\rho^2-u_j^d$ is negative. Choose $L$ with $L^2>\max_j(v_j^d-\rho^2)/g(v_j)^2$, and $L>1$ when $2k=d$. Then every $S_b(v_j)$ is positive, and $S_b$ has a root in each of the $2k$ intervals between consecutive points of the displayed sequence. Its leading coefficient is $L^2$ or $L^2-1$, both positive, so its degree is exactly $2k$ and these are all its roots. They are positive and simple, and the signs of $S_b'$ alternate. This gives $dN$ simple quaternionic roots, half of each local degree.
\end{proof}

For instance, with $d=4$ and $k=2$, the equation
\begin{equation}\label{eq:norm-square-quartic}
q^4+10(|q|^2-2)(|q|^2-4)=i
\end{equation}
has sixteen simple roots, eight of each local degree. Indeed, here $S_b(u)=100(u-2)^2(u-4)^2+1-u^4$ takes the values $900$, $-15$, $20$, $-255$ and $276$ at $u=1,2,3,4,5$. Its leading coefficient $99$ is positive, so its four roots are simple and lie one in each interval between consecutive integers from one to five. The map is conjugation-invariant, has real coordinate degree four, and attains the planar bound $dN=16$ of Proposition~\ref{prop:central}. For even $d$ and $\psi(u)=\phi(u^{d/2})$ one recovers $P_{d,\phi}$, with $S_b(u)=R_b(u^{d/2})$. The positive roots of $S_b$ then correspond to those of $R_b$, which explains the smaller bound of Theorem~\ref{thm:count}.

\section{Persistence and ambient comparisons}\label{sec:stability}
The following elementary coercivity lemma applies to arbitrary real polynomial maps. Its role here is to extend the radial extremal counts to small neighborhoods in the full coefficient space.

\begin{lemma}\label{lem:stability}
Let $F:\RR^4\to\RR^4$ be a polynomial map of degree at most $N$ with finitely many zeros, all simple. Suppose there are $c>0$ and $R>0$ such that
\begin{equation}\label{eq:coercive}
|F(q)|\geq c|q|^N,\qquad |q|\geq R.
\end{equation}
Then every sufficiently small coefficient perturbation among maps of degree at most $N$ has exactly the same number of zeros, all simple and with the same local degrees. The assertion includes a map with no zeros.
\end{lemma}
\begin{proof}
In this finite-dimensional coefficient space there is a constant $C_N$ such that a coefficient perturbation of size $\delta$ changes the map by at most $C_N\delta(1+|q|^N)$. For $|q|\geq\max\{R,1\}$ and sufficiently small $\delta$, inequality~\eqref{eq:coercive} therefore excludes new zeros uniformly.

Around each zero, the implicit function theorem with coefficients as parameters gives a neighborhood in which the perturbed map has a unique simple zero. Shrink these neighborhoods to have disjoint compact closures within the uniqueness neighborhoods. On the remaining compact subset of the ball the original map has positive minimum norm, so sufficiently small coefficient perturbations have no zeros there. Each derivative determinant retains its sign by continuity. The same compact-minimum argument applies when the original zero set is empty.
\end{proof}

\begin{proposition}[Persistence of the radial extremal count]\label{prop:stable}
For each $d=2m\geq2$ and $r\geq1$, the examples attaining the bound in Theorem~\ref{thm:count} can be chosen so that $P_{d,\phi}-b$ has a neighborhood in the coefficient space of all real polynomial maps $\RR^4\to\RR^4$ of degree at most $N=dr$ in which every map has exactly $2N$ simple zeros. Half have local degree $+1$ and half have local degree $-1$.
\end{proposition}
\begin{proof}
Use the sharpness examples $F=P_{d,\phi_L}-I\rho$. For $r\geq2$, the leading term is $L|q|^{dr}$, so~\eqref{eq:coercive} holds with $N=dr$. For $r=1$ and $L>1$, the reverse triangle inequality gives $|F(q)|\geq(L-1)|q|^d-O(1)$. Lemma~\ref{lem:stability} applies in both cases. The local degrees were determined in the sharpness proof.
\end{proof}

The same argument applies to the sharpness examples of Proposition~\ref{prop:norm-square}. Their leading term is $L|q|^{2k}$ when $2k>d$, and $|q^d+L|q|^d|\geq(L-1)|q|^d$ when $2k=d$. Their $dN$ simple roots therefore persist on an open set of polynomial maps of degree at most $N$.

\begin{remark}[Ambient counts]\label{rem:ambient}
The full coefficient space in~\eqref{eq:coefficient-space} admits at most $N^4$ isolated complex common zeros after complexification, counted with multiplicity \cite{fulton}. In particular, it admits at most $N^4$ simple real zeros. A finite complex zero set also suffices for this bound on all distinct real zeros. Real isolation alone does not imply complex isolation.\footnote{For instance, $F:\RR^2\to\RR^2$ defined by $F(x,y)=(x^2+y^2,0)$ has the origin as its only real zero. Its complexified zero set is the union of the lines $x=iy$ and $x=-iy$.} Without complex isolation, Milnor's theorem for a real zero set defined by polynomials of degree at most $N$ bounds the sum of its Betti numbers by $N(2N-1)^3$ \cite[Theorem~2]{milnorbetti}. This also bounds its connected components and isolated points.

The elementary product construction attains the bound. Put $p_N(t)=\prod_{j=1}^N(t-j)$ and, for $q=x_0+x_1i+x_2j+x_3k$, define
\begin{equation}\label{eq:grid-map}
H_N(q)=p_N(x_0)+p_N(x_1)i+p_N(x_2)j+p_N(x_3)k.
\end{equation}
For $N\geq2$, this coordinate grid construction lies outside the one-word class with leading term $q^N$. Its entire complex zero set is $\{1,\ldots,N\}^4$, and its Jacobian on this grid is diagonal with non-zero entries $p_N'(x_\ell)$. The norm inequality $(\sum_{\ell=0}^3x_\ell^{2N})^{1/2}\geq2^{1-N}(\sum_{\ell=0}^3x_\ell^2)^{N/2}$ for its leading term gives
\[
|H_N(q)|\geq2^{1-N}|q|^N-O(|q|^{N-1}+1).
\]
For even $N$, the signs $\sigma(j)=\sgn p_N'(j)=(-1)^{N-j}$ are equally divided, and the signed grid count is $(\sum_{j=1}^N\sigma(j))^4=0$. Thus half of its $N^4$ roots have each local degree. At the other extreme, $F_N(q)=q^N+2|q|^N+1$ has $\Rea F_N(q)\geq1+|q|^N$ and no roots. Lemma~\ref{lem:stability} gives open neighborhoods with these counts. Together with Proposition~\ref{prop:stable}, these elementary comparisons produce the counts $0$, $2N$, and $N^4$ on open subsets for every even $N\geq2$. The examples of Proposition~\ref{prop:norm-square} add the counts $dN$ for $1\leq d\leq N$. They describe local persistence with a fixed degree bound.
\end{remark}

\section{How spherical root sets break into isolated roots}\label{sec:spherical-splitting}
The real fibers described in Subsection~\ref{subsec:real-targets} contain conjugacy two-spheres. Their tangent directions lie in the kernel of the derivative, so the simple-root persistence lemma does not apply. A normal reduction instead leaves a tangent vector field on the sphere. Its zeros determine the isolated roots after perturbation, and its indices determine their local degrees. We first treat a spherical component with an invertible normal derivative, then construct polynomial perturbations realizing every positive even local root count when the perturbation degree is allowed to increase. Section~\ref{sec:eigen} specializes the reduction to imaginary perturbations, obtaining sharp local bounds under a fixed perturbation degree and exact global root counts. The reduction uses the implicit function theorem, while the signed count is an application of the Poincar\'e--Hopf theorem \cite{milnor}.

\subsection{A spherical component and its normal derivative}\label{subsec:spherical-normal}
Fix a real target $a$ and put $F_0(q)=P_{d,\phi}(q)-a$. We use $\mathrm i$ as the imaginary unit of an auxiliary complex plane, identified with $\CC_I=\RR+I\RR\subset\HH$ by $x+\mathrm i y\mapsto x+Iy$ for each unit imaginary quaternion $I$. This planar unit is distinct from the quaternionic basis element $i$, and $\Ima_{\mathrm i}$ denotes its real coefficient. Choose $x_0\in\RR$ and $y_0>0$ such that, for $z_0=x_0+\mathrm i y_0$,
\begin{equation}\label{eq:sphere-data}
z_0^d=\sigma t_0,\qquad t_0=|z_0|^d>0,\qquad \sigma\in\{-1,1\},\qquad a=\sigma t_0+\phi(t_0).
\end{equation}
Then $F_0$ vanishes on
\[
\Sigma=\{x_0+y_0J\mid J\in\mathbb S^2\},\qquad \mathbb S^2=\{J\in\Ima\HH\mid |J|=1\}.
\]
Near this sphere, write $q=x+yJ$ with $y>0$. Define $A,B$ by
\[
A(x,y)=\Rea((x+\mathrm i y)^d)+\phi((x^2+y^2)^{d/2})-a,\qquad B(x,y)=\Ima_{\mathrm i}((x+\mathrm i y)^d).
\]
Thus $F_0(x+yJ)=A(x,y)+JB(x,y)$. Set
\begin{equation}\label{eq:sphere-normal-matrix}
L_\Sigma=D_{(x,y)}(A,B)(x_0,y_0),\qquad \kappa=1+\sigma\phi'(t_0).
\end{equation}
\begin{lemma}[Normal non-degeneracy]\label{lem:sphere-normal}
The normal determinant is
\begin{equation}\label{eq:sphere-normal-det}
\det L_\Sigma=d^2|z_0|^{2d-2}\kappa.
\end{equation}
If $\kappa\neq0$, then $\Sigma$ is an isolated component of $F_0^{-1}(0)$ in a sufficiently small tubular neighborhood, $DF_0$ has rank two along $\Sigma$, and its kernel is exactly $T\Sigma$.
\end{lemma}
\begin{proof}
On the complex slice, the derivative of $z\mapsto z^d$ has determinant $d^2|z_0|^{2d-2}>0$. At the real value $w=\sigma t_0$, the derivative of the planar map $w\mapsto w+\phi(|w|)$ of $\CC$, taken in the real and imaginary coordinates of this auxiliary plane, is $\operatorname{diag}(1+\sigma\phi'(t_0),1)$. Their product gives~\eqref{eq:sphere-normal-det}. The tangent derivative vanishes because $A(x_0,y_0)=B(x_0,y_0)=0$. Invertibility of $L_\Sigma$ proves the rank and kernel assertions. The inverse function theorem isolates $(x_0,y_0)$ as a zero of $(A,B)$, and therefore isolates the corresponding sphere.
\end{proof}

\subsection{Reduction to a tangent vector field}\label{subsec:spherical-reduction}
For $J\in\mathbb S^2$, let $\Pi_Jv=v-\langle v,J\rangle J$ be orthogonal projection onto $T_J\mathbb S^2$. Let $H$ be a smooth quaternion-valued map on a neighborhood of $\Sigma$ and consider $F_\varepsilon=F_0+\varepsilon H$. Define
\begin{equation}\label{eq:sphere-field}
V_H(J)=\Pi_J\Ima H(x_0+y_0J).
\end{equation}
A zero of this field is non-degenerate when its derivative as a tangent section is invertible. Its index is then $+1$ or $-1$, according to the sign of that determinant.

\begin{theorem}[Splitting and the signed count]\label{thm:sphere-splitting}
Fix $H$ and assume $\kappa\neq0$ and that every zero of $V_H$ is non-degenerate. For this $H$, choose a sufficiently small bounded isolating tubular neighborhood $U$ of $\Sigma$, independent of $\varepsilon$. There is then $\varepsilon_0(H,U)>0$ such that, for $0<|\varepsilon|<\varepsilon_0(H,U)$, the roots of $F_\varepsilon$ in $U$ are simple and correspond bijectively to the zeros of $V_H$. The root corresponding to $J_*$ has local degree
\begin{equation}\label{eq:sphere-local-sign}
\operatorname{sgn}(\kappa)\operatorname{ind}_{J_*}(V_H).
\end{equation}
In particular, their signed sum is $2\operatorname{sgn}(\kappa)$. Their number is a positive even integer. If this number is $2\ell+2$, then the counts $(n_+,n_-)$ of positive and negative local degrees are $(\ell+2,\ell)$ for $\kappa>0$ and $(\ell,\ell+2)$ for $\kappa<0$.
\end{theorem}
\begin{proof}
Write $h_0(x,y,J)=\Rea H(x+yJ)$ and $h_\parallel(x,y,J)=\langle\Ima H(x+yJ),J\rangle$. The real and $J$-parallel equations are
\begin{equation}\label{eq:sphere-normal-equations}
A(x,y)+\varepsilon h_0(x,y,J)=0,\qquad B(x,y)+\varepsilon h_\parallel(x,y,J)=0.
\end{equation}
The implicit function theorem and compactness of $\mathbb S^2$ give unique smooth solutions $x_\varepsilon(J),y_\varepsilon(J)$ close to $x_0,y_0$, uniformly in $J$. The remaining equation, after division by $\varepsilon\neq0$, is the tangent-section equation
\begin{equation}\label{eq:sphere-reduced-section}
V_\varepsilon(J)=\Pi_J\Ima H(x_\varepsilon(J)+y_\varepsilon(J)J)=0.
\end{equation}
This section converges to $V_H$ in $C^1$. Each non-degenerate zero therefore persists uniquely. Away from disjoint neighborhoods of these zeros, $|V_H|$ has a positive minimum, excluding further zeros. A smaller tube and a compact-annulus argument exclude roots outside the normal-coordinate neighborhood used for the reduction.

For the sign, choose an orthonormal tangent frame $\tau_1,\tau_2$ with $(J,\tau_1,\tau_2)$ positively oriented. Use $(x,y)$ and the tangent coordinates in the domain, and the frame $(1,J,\tau_1,\tau_2)$ in the target. At a root, write the derivative of the real, sphere-parallel and tangent equations as $\begin{pmatrix}M&C\\R&S\end{pmatrix}$, where $M$ is the invertible normal block. Differentiating~\eqref{eq:sphere-normal-equations} along the graph $(x_\varepsilon(J),y_\varepsilon(J))$ gives the normal-coordinate derivative $-M^{-1}C$. Hence the derivative of the remaining tangent equation is $S-RM^{-1}C=\varepsilon DV_\varepsilon(J)$. Terms from differentiating the moving target frame vanish because $F_\varepsilon(q)=0$. The target frame has the orientation of $(1,i,j,k)$. The parametrization sends the ordered domain frame to $(1,J,y\tau_1,y\tau_2)$ and has determinant $y^2>0$ in this target frame. Thus
\begin{equation}\label{eq:sphere-schur-determinant}
\det DF_\varepsilon(q)
=\frac{\det M}{y_\varepsilon(J)^2}\det(S-RM^{-1}C)
=\frac{\varepsilon^2\det M}{y_\varepsilon(J)^2}\det DV_\varepsilon(J).
\end{equation}
The factor $\varepsilon^2$ comes from the two tangent directions and is positive for either sign of $\varepsilon$. The normal determinant retains the sign of $\det L_\Sigma$, and the tangent determinant retains the index of the corresponding zero of $V_H$. This proves~\eqref{eq:sphere-local-sign}. The Poincar\'e--Hopf theorem gives
\[
\sum_{V_H(J)=0}\operatorname{ind}_J(V_H)=\chi(\mathbb S^2)=2.
\]
The assertions about the counts follow because every index is $+1$ or $-1$.
\end{proof}

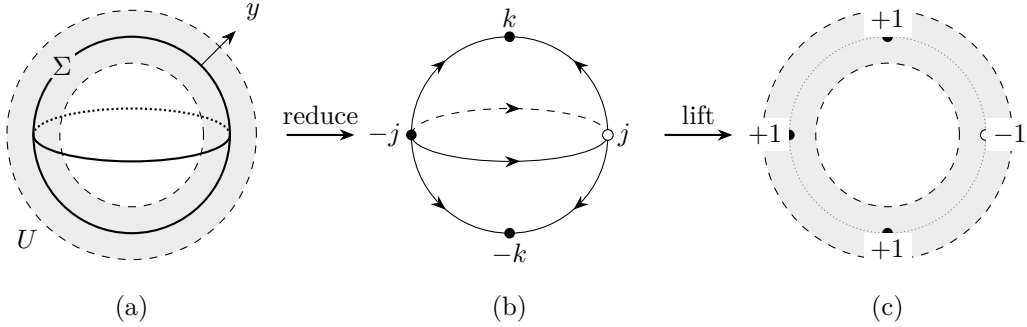
\begin{figure}[htbp]
\centering
\begin{tikzpicture}[>={Stealth[length=2mm]},font=\small]
\begin{scope}
\fill[gray!15,even odd rule] (0,0) circle (1.65) (0,0) circle (0.95);
\draw[dashed] (0,0) circle (1.65);
\draw[dashed] (0,0) circle (0.95);
\draw[thick] (0,0) circle (1.3);
\draw[thick] (-1.3,0) arc (180:360:1.3 and 0.35);
\draw[thick,densely dotted] (1.3,0) arc (0:180:1.3 and 0.35);
\draw[->] (45:1.3) -- (45:1.95) node[above right] {$y$};
\node[fill=gray!15,inner sep=1pt] at (135:1.3) {$\Sigma$};
\node at (225:1.95) {$U$};
\node at (0,-2.3) {(a)};
\end{scope}
\draw[->,thick] (2.05,0) -- (2.95,0) node[midway,above] {reduce};
\begin{scope}[xshift=5cm]
\draw[midarrow] (-1.3,0) arc (180:90:1.3);
\draw[midarrow] (-1.3,0) arc (180:270:1.3);
\draw[midarrow] (1.3,0) arc (0:90:1.3);
\draw[midarrow] (1.3,0) arc (0:-90:1.3);
\draw[midarrow] (-1.3,0) arc (180:360:1.3 and 0.35);
\draw[dashed,midarrow] (-1.3,0) arc (180:0:1.3 and 0.35);
\fill (0,1.3) circle (2pt) node[above] {$k$};
\fill (0,-1.3) circle (2pt) node[below] {$-k$};
\fill (-1.3,0) circle (2pt) node[left] {$-j$};
\draw[fill=white] (1.3,0) circle (2pt) node[right] {$j$};
\node at (0,-2.3) {(b)};
\end{scope}
\draw[->,thick] (7.05,0) -- (7.95,0) node[midway,above] {lift};
\begin{scope}[xshift=10cm]
\fill[gray!15,even odd rule] (0,0) circle (1.65) (0,0) circle (0.95);
\draw[dashed] (0,0) circle (1.65);
\draw[dashed] (0,0) circle (0.95);
\draw[gray,densely dotted] (0,0) circle (1.3);
\fill (0,1.3) circle (2pt) node[above,fill=white,inner sep=2.5pt] {$+1$};
\fill (0,-1.3) circle (2pt) node[below,fill=white,inner sep=2.5pt] {$+1$};
\fill (-1.3,0) circle (2pt) node[left,fill=white,inner sep=2.5pt] {$+1$};
\draw[fill=white] (1.3,0) circle (2pt) node[right,fill=white,inner sep=2.5pt] {$-1$};
\node at (0,-2.3) {(c)};
\end{scope}
\end{tikzpicture}
\caption{Normal elimination and spherical splitting for $\kappa>0$, with the first field $V_1$ of~\eqref{eq:sphere-model-field}, $\zeta_1=0$. Panel (a) shows a section of the isolating tube $U$ around $\Sigma$, with the real normal direction suppressed. In (b), the sinks $\pm k$ and the source $-j$ have index $+1$, while the saddle $j$ has index $-1$. In (c), each zero lifts to one simple root whose local degree equals its index. The signed sum is two. The roots are drawn on $\Sigma$, as for $H_1$ in~\eqref{eq:sphere-polynomial-perturbation}. General perturbations with this first field move them within $U$.}
\label{fig:sphere-splitting}
\end{figure}

\begin{corollary}[A local degree obstruction to removing the sphere]\label{cor:sphere-degree}
Under the normal non-degeneracy hypothesis, the local Brouwer degree on a sufficiently small tube is
\[
\deg(F_0,U,0)=2\operatorname{sgn}(\kappa).
\]
Every continuous perturbation sufficiently small on $\partial U$ has a zero in $U$. If all its zeros in $U$ are simple, then their signed sum is $2\operatorname{sgn}(\kappa)$ and there are at least two.
\end{corollary}
\begin{proof}
Choose a constant unit imaginary perturbation $H=e$. Its tangent field $e-\langle e,J\rangle J$ has exactly two non-degenerate zeros, at $J=\pm e$, each of index $+1$. Theorem~\ref{thm:sphere-splitting} computes the degree of $F_0+\varepsilon e$. The straight-line homotopy to $F_0$ avoids zero on $\partial U$ for small $\varepsilon$. The same boundary estimate applies to any perturbation whose norm on $\partial U$ is less than $\min_{\partial U}|F_0|$. Homotopy invariance and the regular-value formula prove the result.
\end{proof}

\begin{remark}[Genericity and the next reduced term]\label{rem:sphere-genericity}
In any full finite-dimensional space of polynomial perturbations containing the constants, the fields $V_H$ have only non-degenerate zeros for an open dense set of coefficients. Indeed, variation of the imaginary constant term maps onto $T_J\mathbb S^2$ at every $J$. Parametric transversality therefore applies to the universal section, and compactness makes non-degeneracy open \cite[Chapter~II, \S4]{golubitskyguillemin}. This assertion concerns generic perturbations near the chosen sphere.

The reduction also gives an expansion when the first field is degenerate. Write $q_J=x_0+y_0J$ and let $L_\Sigma$ be the constant matrix from~\eqref{eq:sphere-normal-matrix}, independent of $J$. Define real functions $\xi_x,\xi_y$ by
\begin{equation}\label{eq:sphere-normal-displacement}
\begin{pmatrix}\xi_x(J)\\\xi_y(J)\end{pmatrix}
=-L_\Sigma^{-1}\begin{pmatrix}\Rea H(q_J)\\\langle\Ima H(q_J),J\rangle\end{pmatrix}.
\end{equation}
The dependence of $\xi_x,\xi_y$ on $J$ enters solely through the two displayed components of $H(q_J)$. Then, uniformly on the sphere,
\begin{equation}\label{eq:sphere-second-section}
V_\varepsilon(J)=V_H(J)+\varepsilon\Pi_J\Ima\bigl(DH(q_J)[\xi_x(J)+\xi_y(J)J]\bigr)+O(\varepsilon^2).
\end{equation}
The expansion holds for every smooth $H$. If $V_H$ vanishes identically and the coefficient of $\varepsilon$ in~\eqref{eq:sphere-second-section} has only non-degenerate zeros, then dividing the full tangent equation by $\varepsilon^2$ gives the same correspondence with that second field and the same signed-count law. Near an isolated degenerate zero of $V_H$, the dependence on both tangent and normal coordinates requires a separate local analysis.

For a conjugation-equivariant perturbation, $H(x+yJ)$ lies in $\RR+J\RR$, so the exact section $V_\varepsilon$ and every reduced term vanish identically. The two normal equations are independent of $J$ and produce a nearby sphere. Thus the isolated-root constructions below require perturbations that break this symmetry.
\end{remark}

\subsection{Polynomial realization of every even local count}\label{subsec:spherical-realization}
We identify $\Ima\HH$ with oriented Euclidean three-space and write $e_1=i$, $e_2=j$, $e_3=k$. Fix an integer $\ell\geq0$ and distinct numbers $-1<\zeta_1<\cdots<\zeta_\ell<1$. Put $p_\ell(z)=\prod_{j=1}^{\ell}(z-\zeta_j)$, with $p_0=1$. The tangent polynomial field
\begin{equation}\label{eq:sphere-model-field}
V_\ell(J)=p_\ell(J_3)(e_3-J_3J)+J_1(e_3\times J)
\end{equation}
has two polar zeros and one opposite-index pair on each latitude $J_3=\zeta_j$.

\begin{proposition}[Complete realization with variable perturbation degree]\label{prop:sphere-realization}
For every normally non-degenerate spherical component $\Sigma$ as above and every $\ell\geq0$, there is a real polynomial perturbation $H_\ell$ of degree at most $\ell+2$ such that $F_0+\varepsilon H_\ell$ has exactly $2\ell+2$ simple roots in a fixed tube around $\Sigma$ for every $\varepsilon\neq0$. Their signed counts are those in Theorem~\ref{thm:sphere-splitting}. Consequently every positive even local count is attainable when no common bound on perturbation degree is imposed.
\end{proposition}
\begin{proof}
For $v=\Ima q$, define the purely imaginary polynomial
\begin{equation}\label{eq:sphere-polynomial-perturbation}
H_\ell(q)=\frac{p_\ell(v_3/y_0)}{y_0^2}\bigl(|v|^2e_3-v_3v\bigr)
+\frac{v_1}{y_0^2}(e_3\times v).
\end{equation}
It has degree at most $\ell+2$, is perpendicular to $v$ everywhere, and restricts to $V_\ell(J)$ on $v=y_0J$. Thus its real and $J$-parallel components vanish identically in the tube. The normal equations remain $A=B=0$, whose only nearby solution is $(x_0,y_0)$. The exact remaining equation is $\varepsilon V_\ell(J)=0$.

Away from the poles, use coordinates $z=J_3$ and $\theta$ with $J_1=\sqrt{1-z^2}\cos\theta$. The field becomes
\[
V_\ell=p_\ell(z)(1-z^2)\frac{\partial}{\partial z}
+\sqrt{1-z^2}\cos\theta\frac{\partial}{\partial\theta}.
\]
Its zeros have $z=\zeta_j$ and $\theta=\pi/2,3\pi/2$. The two diagonal derivatives are $p_\ell'(\zeta_j)(1-\zeta_j^2)$ and $-\sqrt{1-\zeta_j^2}\sin\theta$, so the two zeros on each latitude are non-degenerate and have opposite indices. At a zero, changing to an orthonormal tangent frame conjugates the derivative and leaves its determinant unchanged. Near either pole, the second term in~\eqref{eq:sphere-model-field} vanishes to second order. The first has linear part $-p_\ell(1)\operatorname{Id}$ at the north pole and $p_\ell(-1)\operatorname{Id}$ at the south pole. Each pole is non-degenerate with index $+1$. There are therefore $\ell+2$ positive indices and $\ell$ negative indices. The determinant argument in Theorem~\ref{thm:sphere-splitting} applies for every $\varepsilon\neq0$ here, since the normal equations are unchanged. This proves the proposition.
\end{proof}

\subsection{An exact global example based on the quadratic power}\label{subsec:spherical-global}
The zero set of $q^2+1$ is the unit imaginary sphere. Here the radial term is $\phi=0$, for which the preceding normal reduction remains valid and gives $\kappa=1$.

\begin{theorem}[A proper polynomial with any positive even root count]\label{thm:sphere-global}
Use~\eqref{eq:sphere-polynomial-perturbation} with $y_0=1$ and define
\[
G_{\varepsilon,\ell}(q)=q^2+1+\varepsilon H_\ell(q),\qquad \varepsilon\neq0.
\]
This map has real coordinate degree $\ell+2$, is proper of Brouwer degree two, and has exactly $2\ell+2$ roots, all simple. The roots are
\begin{equation}\label{eq:sphere-explicit-roots}
q=\pm k
\qquad\text{and}\qquad
q=\pm\sqrt{1-\zeta_j^2}\,j+\zeta_j k
\quad(1\leq j\leq\ell).
\end{equation}
Exactly $\ell+2$ roots have local degree $+1$ and $\ell$ have local degree $-1$.
\end{theorem}
\begin{proof}
Write $q=x+v$ with $x\in\RR$ and $v\in\Ima\HH$. Then
\[
G_{\varepsilon,\ell}(q)=x^2-|v|^2+1+2xv+\varepsilon H_\ell(v).
\]
Taking the inner product of its imaginary part with $v$ gives $2x|v|^2=0$. A root cannot have $v=0$, since $x^2+1>0$. Hence $x=0$, and the real equation gives $|v|=1$. The remaining equation is $V_\ell(v)=0$, which gives exactly~\eqref{eq:sphere-explicit-roots}. The normal matrix of $(x^2-y^2+1,2xy)$ at $(0,1)$ has determinant four, so the index calculation in Proposition~\ref{prop:sphere-realization} proves simplicity and the stated signs.

Orthogonality of $H_\ell(v)$ and $v$ also gives
\begin{align*}
|G_{\varepsilon,\ell}(q)|^2
&=(x^2-|v|^2+1)^2+4x^2|v|^2+\varepsilon^2|H_\ell(v)|^2\\
&=|q^2+1|^2+\varepsilon^2|H_\ell(v)|^2
\geq (|q|^2-1)^2.
\end{align*}
This proves properness uniformly along the homotopy replacing $\varepsilon$ by $s\varepsilon$ for $0\leq s\leq1$. The degree is therefore that of $q^2+1$, namely two. For $\ell>0$, the non-zero leading term of $H_\ell$ gives coordinate degree $\ell+2$. For $\ell=0$, the real part is quadratic.
\end{proof}

For instance, $\ell=1$ and $\zeta_1=0$ give a cubic perturbation with exactly four roots $\pm j,\pm k$, three positive and one negative, as in Figure~\ref{fig:sphere-splitting}. As $\varepsilon$ tends to zero, these maps converge to $q^2+1$ on every compact set. The growing root counts require increasing polynomial degree. For $\ell>0$, the quadratic growth estimate does not give the degree-$\ell+2$ coercivity of Lemma~\ref{lem:stability}, and the exact count is unstable in the full coefficient space. To see this, put
\begin{equation}\label{eq:sphere-global-escaping-root}
\widetilde G_{\delta}(q)=G_{\varepsilon,\ell}(q)-\delta(\Rea q)^{\ell+2},\qquad \delta>0.
\end{equation}

Fix $\ell>0$ and $\varepsilon\neq0$. Orthogonality of $H_\ell(v)$ and $v$ gives $\langle\Ima\widetilde G_\delta(x+v),v\rangle=2x|v|^2$. Every non-real root therefore has $x=0$, so the non-real roots are precisely~\eqref{eq:sphere-explicit-roots}, with their derivatives unchanged. On the real axis the equation is $x^2+1=\delta x^{\ell+2}$. The ratio $x^{-\ell}+x^{-\ell-2}$ decreases strictly from infinity to zero on $x>0$, giving a unique positive root $x_\delta$, which tends to infinity as $\delta$ decreases to zero. For even $\ell$ there is also the root $-x_\delta$, while for odd $\ell$ there is no negative real root. Since $DH_\ell(0)=0$, the determinant at either real root is $(-\ell x-(\ell+2)/x)(2x)^3<0$. Thus every root is simple, and each additional real root has local degree $-1$.

The estimate $|\widetilde G_\delta(x+v)|\geq2|x||v|$ also proves properness. A bounded image along an unbounded sequence would force $|x||v|$ to remain bounded. Along a subsequence, either $|x|\to\infty$ and $|v|\to0$, or $|v|\to\infty$ and $x\to0$. In both cases the real part $x^2-|v|^2+1-\delta x^{\ell+2}$ is unbounded. The complete count and the regular-value formula therefore give
\begin{equation}\label{eq:sphere-global-perturbed-count}
\bigl(\#\widetilde G_\delta^{-1}(0),\,\deg\widetilde G_\delta\bigr)=
\begin{cases}
(2\ell+3,\,1), & \ell\text{ odd},\\
(2\ell+4,\,0), & \ell\text{ even}.
\end{cases}
\end{equation}
Arbitrarily small coefficient perturbations thus change the complete root count through roots arriving from infinity.

\section{Tensor eigenvectors and polynomials with one leading word}\label{sec:eigen}
When the perturbation depends only on the imaginary part of $q$ and takes imaginary values, the reduction of Section~\ref{sec:spherical-splitting} becomes exact and global. We use it here to obtain sharp local bounds under a fixed perturbation degree, a formula for the Brouwer degree of quadratic families, and lower bounds for the number of roots of polynomials with a single leading word.

\subsection{Eigenvectors of imaginary perturbations}\label{subsec:spherical-eigen}
The field in~\eqref{eq:sphere-field} takes a particularly transparent form when the perturbation depends only on the imaginary part of $q$ and takes imaginary values. In this case the whole root problem separates, in an exact and global manner, into an eigenvector problem on $\Ima\HH$ and a planar problem on a complex slice. We find this separation illuminating also for the last observation of Remark~\ref{rem:sphere-genericity}, since a perturbation respecting conjugation admits every direction as an eigenvector.

Let $h\colon\Ima\HH\to\Ima\HH$ be a polynomial map, homogeneous of degree $D\geq1$. A unit eigenvector of $h$ is a point $J\in\mathbb S^2$ with $h(J)=\mu J$ for some real number $\mu$, its eigenvalue. We write $\mathcal E(h)$ for the set of unit eigenvectors and put
\[
V_h(J)=\Pi_Jh(J),\qquad \|h\|=\max_{J\in\mathbb S^2}|h(J)|.
\]
Thus $\mathcal E(h)$ is the zero set of the tangent field $V_h$. At a zero $J$ with eigenvalue $\mu$, differentiation of $V_h(J')=h(J')-\langle h(J'),J'\rangle J'$ along a tangent vector gives the intrinsic derivative
\begin{equation}\label{eq:eigen-derivative}
K_J=\bigl(\Pi_JDh(J)-\mu\bigr)\big|_{T_J\mathbb S^2}.
\end{equation}
The real coordinates are two-sided expressions by~\eqref{eq:coefficient-space}, so $q\mapsto h(\Ima q)$ is a two-sided polynomial of degree $D$. Non-degeneracy of $V_h$ is open and dense among homogeneous maps of degree $D$ and among gradients of homogeneous forms of degree $D+1$. At each $J$, the variation $\delta h(v)=\langle J,v\rangle^D c$, with $c\in T_J\mathbb S^2$, gives $\delta V_h(J)=c$. In the gradient class $h=\nabla W$, the variation $\delta W(v)=\langle J,v\rangle^D\langle c,v\rangle$ gives the same tangent value. Parametric transversality and compactness give the claim, as in Remark~\ref{rem:sphere-genericity}.

\begin{lemma}[Eigenvector reduction]\label{lem:eigen-reduction}
Let $F\colon\HH\to\HH$ be a polynomial map, and suppose there are smooth real functions $A,B$ on $\RR\times(0,\infty)$ with $F(x+yJ)=A(x,y)+JB(x,y)$ for all $x\in\RR$, $y>0$ and $J\in\mathbb S^2$. Let $h$ be as above and put $G(q)=F(q)+h(\Ima q)$. The real roots of $G$ coincide with those of $F$. A point $q=x+yJ$ with $y>0$ is a root of $G$ exactly when $J\in\mathcal E(h)$ and, with $\mu$ the eigenvalue of $J$,
\begin{equation}\label{eq:eigen-planar}
A(x,y)=0,\qquad B(x,y)+\mu y^D=0.
\end{equation}
At such a root,
\begin{equation}\label{eq:eigen-det}
\det DG(q)=y^{2D-2}\,\det D_{(x,y)}\bigl(A,B+\mu y^D\bigr)(x,y)\,\det K_J.
\end{equation}
The root is simple exactly when both determinants on the right side are non-zero. Its local degree is then
\[
\sgn\det D_{(x,y)}\bigl(A,B+\mu y^D\bigr)(x,y)\cdot\operatorname{ind}_J(V_h).
\]
\end{lemma}
\begin{proof}
The map $h$ vanishes at the origin, being homogeneous of positive degree, so $G$ and $F$ agree on $\RR$. For $y>0$, homogeneity gives $h(yJ)=y^Dh(J)$ and hence
\[
G(x+yJ)=A(x,y)+JB(x,y)+y^Dh(J).
\]
The real part of this quaternion is $A(x,y)$. Its imaginary part vanishes exactly when $h(J)=-y^{-D}B(x,y)J$, that is, when $J$ is a unit eigenvector whose eigenvalue satisfies $B(x,y)+\mu y^D=0$.

For the derivative, choose an orthonormal basis $\tau_1,\tau_2$ of $T_J\mathbb S^2$. The map $(x,y,J')\mapsto x+yJ'$ is a local diffeomorphism near $(x,y,J)$, and it sends the four coordinate directions to $1$, $J$, $y\tau_1$ and $y\tau_2$. We express their images under $DG(q)$ in the orthonormal basis $(1,J,\tau_1,\tau_2)$ of $\HH$. The real component of $G$ depends only on $(x,y)$. The derivative of $JB+y^Dh(J)$ in $y$ equals $J(B_y+D\mu y^{D-1})$, and the derivatives in $x$ and $y$ have no tangent component, because $h(J)$ is parallel to $J$. Along $y\tau_a$, the derivative equals $B\tau_a+y^DDh(J)\tau_a$. Since $B=-\mu y^D$ at the root, its tangent component equals $y^DK_J\tau_a$ by~\eqref{eq:eigen-derivative}. The resulting matrix is block upper triangular,
\[
\begin{pmatrix}
D_{(x,y)}(A,B+\mu y^D) & \ast\\
0 & y^DK_J
\end{pmatrix}.
\]
The first row of the upper-right block vanishes, since the real component of $G$ depends only on $(x,y)$. Its second row consists of the entries $y^D\langle Dh(J)\tau_a,J\rangle$, which need not vanish when $h$ is not a gradient. These entries play no role in the determinant, because the lower-left block vanishes identically. The passage to the orthonormal basis in the domain divides the last two columns by $y$, so the tangent block becomes $y^{D-1}K_J$. Its determinant equals $(y^{D-1})^2\det K_J=y^{2D-2}\det K_J$, since $T_J\mathbb S^2$ has dimension two. We now use the same orthonormal basis $(1,J,\tau_1,\tau_2)$ in the domain and codomain, so its orientation relative to $(1,i,j,k)$ cancels in the determinant. This proves~\eqref{eq:eigen-det}. The index of a non-degenerate zero of $V_h$ is the sign of $\det K_J$.
\end{proof}

For $H(q)=h(\Ima q)$, the field~\eqref{eq:sphere-field} equals $y_0^DV_h$. Lemma~\ref{lem:eigen-reduction} holds for every size of the perturbation, whereas Theorem~\ref{thm:sphere-splitting} describes small perturbations. In the small limit the planar determinant in~\eqref{eq:eigen-det} tends to $\det L_\Sigma$, whose sign is $\sgn\kappa$, and one recovers~\eqref{eq:sphere-local-sign}.

\begin{lemma}[Antipodal pairs and parity]\label{lem:eigen-antipodal}
Let $h$ be homogeneous of degree $D\geq1$. Then $-\mathcal E(h)=\mathcal E(h)$, and the eigenvalue of $-J$ equals $(-1)^{D+1}$ times that of $J$. A zero $J$ of $V_h$ is non-degenerate exactly when $-J$ is, and in that case $\operatorname{ind}_{-J}(V_h)=\operatorname{ind}_J(V_h)$. Assume that every zero of $V_h$ is non-degenerate. Then $\mathcal E(h)$ consists of $2e(h)$ points, where the number $e(h)$ of real eigenvector lines is odd, and the indices of one representative from each line add up to one.
\end{lemma}
\begin{proof}
Homogeneity gives $h(-J)=(-1)^Dh(J)$, and $\Pi_{-J}=\Pi_J$. Hence $V_h(-J)=(-1)^DV_h(J)$, which proves the first assertion. The antipodal map $a(J)=-J$ has derivative $-\mathrm{Id}$ from $T_J\mathbb S^2$ to $T_{-J}\mathbb S^2=T_J\mathbb S^2$, so the push-forward of $V_h$ by $a$ equals $(-1)^{D+1}V_h$. The index of a non-degenerate zero is invariant under push-forward by a diffeomorphism, and on a surface it is unchanged when the field is multiplied by $-1$. This gives the equality of indices. Under the non-degeneracy assumption, $\mathcal E(h)$ is finite, and the Poincar\'e--Hopf theorem gives $2\sum\operatorname{ind}=\chi(\mathbb S^2)=2$, the sum running over one representative of each line. A sum of $e(h)$ terms equal to $\pm1$ can equal one only when $e(h)$ is odd.
\end{proof}

Under homogeneous imaginary perturbations, the local counts of Theorem~\ref{thm:sphere-splitting} are therefore congruent to two modulo four. The perturbations $H_\ell$ of Proposition~\ref{prop:sphere-realization} are imaginary and mix several degrees, and they realize the counts divisible by four as well.

\subsection{Exact examples of low perturbation degree}\label{subsec:spherical-low-degree}
We return to the unit sphere of $q^2+1$, for which $A(x,y)=x^2-y^2+1$ and $B(x,y)=2xy$. With $h$ replaced by $\varepsilon h$, the planar system~\eqref{eq:eigen-planar} over an eigenvector with eigenvalue $\mu$ becomes
\begin{equation}\label{eq:eigen-quadratic-planar}
x=-\tfrac12\varepsilon\mu y^{D-1},\qquad x^2-y^2+1=0.
\end{equation}
At a solution, the planar determinant in~\eqref{eq:eigen-det} equals
\begin{equation}\label{eq:eigen-quadratic-determinant}
4y^2-(D-1)\,\varepsilon^2\mu^2y^{2D-2}.
\end{equation}

\begin{proposition}[Six roots from a linear perturbation]\label{prop:sphere-linear}
Let $M$ be a real linear map of $\Ima\HH$ with three distinct real eigenvalues $\mu_1<\mu_2<\mu_3$ and unit eigenvectors $u_1,u_2,u_3$. For every real $\varepsilon\neq0$, the map
\[
G_\varepsilon(q)=q^2+1+\varepsilon M(\Ima q)
\]
is proper of Brouwer degree two, and its roots are exactly the six quaternions
\begin{equation}\label{eq:sphere-linear-roots}
q=-\frac{\varepsilon\mu_k}{2}\pm\sqrt{1+\frac{\varepsilon^2\mu_k^2}{4}}\;u_k,\qquad k=1,2,3.
\end{equation}
The two signs correspond to $J=\pm u_k$, with the same positive imaginary radius $y$. All six are simple. The four roots with $k\in\{1,3\}$ have local degree $+1$, and the two roots with $k=2$ have local degree $-1$.
\end{proposition}
\begin{proof}
Here $D=1$ and $\mathcal E(M)=\{\pm u_1,\pm u_2,\pm u_3\}$, the eigenvalue of $\pm u_k$ being $\mu_k$ for either sign. The system~\eqref{eq:eigen-quadratic-planar} has the unique solution $x=-\varepsilon\mu_k/2$ and $y=(1+\varepsilon^2\mu_k^2/4)^{1/2}$, with planar determinant $4y^2>0$. There is no real root, since $x^2+1>0$. The operator $M-\mu_k$ annihilates $u_k$ and induces on the quotient $\Ima\HH/\RR u_k$ an endomorphism with eigenvalues $\mu_j-\mu_k$, for $j\neq k$. The projection $\Pi_{u_k}$ identifies this quotient with $T_{u_k}\mathbb S^2$, and~\eqref{eq:eigen-derivative} shows that $K_{\pm u_k}$ corresponds to the induced endomorphism. Hence
\[
\det K_{\pm u_k}=\prod_{j\neq k}(\mu_j-\mu_k),
\]
which is positive for $k=1,3$ and negative for $k=2$. Lemma~\ref{lem:eigen-reduction} gives simplicity and the local degrees. The estimate $|G_\varepsilon(q)|\geq|q|^2-|\varepsilon|\,\|M\|\,|q|-1$ proves properness, and the Brouwer degree equals the signed count $4-2$.
\end{proof}

The perturbation in Proposition~\ref{prop:sphere-linear} has degree one, whereas Proposition~\ref{prop:sphere-realization} uses degree $\ell+2$ to produce $2\ell+2$ roots. The map $G_\varepsilon$ also satisfies the coercivity hypothesis~\eqref{eq:coercive} with $N=2$, so Lemma~\ref{lem:stability} lets its six simple roots persist on an open set of polynomial maps of degree at most two. For symmetric $M$, the field $V_M$ is the spherical gradient of $\frac12\langle MJ,J\rangle$, and the three pairs of roots lie over its maximum, its saddle and its minimum.

The linear case extends to quadratic perturbations of $q^2+1$, and the root count then depends on the size of the perturbation. The next proposition describes it completely.

\begin{proposition}[Quadratic imaginary perturbations of $q^2+1$]\label{prop:quad-wall}
Let $h\colon\Ima\HH\to\Ima\HH$ be homogeneous of degree two, with every zero of $V_h$ non-degenerate. For $J\in\mathcal E(h)$ write $\mu_J$ for its eigenvalue. For real $\varepsilon\neq0$ put
\[
G_\varepsilon(q)=q^2+1+\varepsilon h(\Ima q).
\]
The roots of $G_\varepsilon$ are the quaternions
\begin{equation}\label{eq:quad-wall-roots}
q_J=\Bigl(1-\frac{\varepsilon^2\mu_J^2}{4}\Bigr)^{-1/2}\Bigl(-\frac{\varepsilon\mu_J}{2}+J\Bigr),\qquad J\in\mathcal E(h),\qquad |\varepsilon\mu_J|<2.
\end{equation}
Each of them is simple, and $q_J$ has local degree $\operatorname{ind}_J(V_h)$. The map is proper exactly when $|\varepsilon\mu_J|\neq2$ for every $J\in\mathcal E(h)$, and its Brouwer degree is then
\begin{equation}\label{eq:quad-wall-degree}
\deg G_\varepsilon=\sum_{J\in\mathcal E(h),\ |\varepsilon\mu_J|<2}\operatorname{ind}_J(V_h).
\end{equation}
\end{proposition}
\begin{proof}
There is no real root, since $h(0)=0$ and $x^2+1>0$. The unit eigenvectors of $\varepsilon h$ are those of $h$, with eigenvalues $\varepsilon\mu_J$, and the field $V_{\varepsilon h}=\varepsilon V_h$ has the same zeros and the same indices as $V_h$, the tangent plane being two-dimensional. By Lemma~\ref{lem:eigen-reduction} and~\eqref{eq:eigen-quadratic-planar} with $D=2$, the roots over $J$ are the points $x+yJ$ with $y>0$, $x=-\frac12\varepsilon\mu_Jy$ and $x^2-y^2+1=0$. The last equation becomes $y^2(1-\varepsilon^2\mu_J^2/4)=1$. It has a solution with $y>0$ exactly when $|\varepsilon\mu_J|<2$, and this solution is unique. This gives~\eqref{eq:quad-wall-roots}, and distinct eigenvectors give distinct roots. The planar determinant equals $y^2(4-\varepsilon^2\mu_J^2)>0$, so Lemma~\ref{lem:eigen-reduction} gives simplicity and the local degrees.

Let $Q_\varepsilon(q)=q^2+\varepsilon h(\Ima q)$, so that $G_\varepsilon=Q_\varepsilon+1$ with $Q_\varepsilon$ homogeneous of degree two. Suppose that $Q_\varepsilon(x+v)=0$ with $x+v\neq0$. Then $x^2=|v|^2$ and $2xv+\varepsilon h(v)=0$. Hence $v\neq0$, and writing $v=|v|J$ one finds $\varepsilon|v|h(J)=-2xJ$. Thus $J\in\mathcal E(h)$, with $\varepsilon\mu_J|v|=-2x$ and therefore $|\varepsilon\mu_J|=2$. Conversely, suppose that $|\varepsilon\mu_J|=2$ for some $J\in\mathcal E(h)$. Then the quaternion $q_*=-\frac12\varepsilon\mu_J+J$ satisfies $Q_\varepsilon(q_*)=\frac14\varepsilon^2\mu_J^2-1-\varepsilon\mu_JJ+\varepsilon\mu_JJ=0$, so $G_\varepsilon(tq_*)=1$ for every real $t$ and properness fails. In the remaining case $Q_\varepsilon$ has a positive minimum $c$ on the unit sphere, and $|G_\varepsilon(q)|\geq c|q|^2-1$. The map is then proper, zero is a regular value, and the regular-value formula gives~\eqref{eq:quad-wall-degree}.
\end{proof}

\begin{remark}[The degree as a step function]\label{rem:quad-wall}
Formula~\eqref{eq:quad-wall-degree} presents the Brouwer degree as a step function of $|\varepsilon|$. It changes only at the values $2/|\mu_J|$, where the roots over the corresponding eigenvectors escape to infinity together. Indeed, the root~\eqref{eq:quad-wall-roots} has norm $\bigl((4+\varepsilon^2\mu_J^2)/(4-\varepsilon^2\mu_J^2)\bigr)^{1/2}$, and as $|\varepsilon\mu_J|$ increases to $2$ it escapes along the ray through the non-trivial zero $q_*=-\frac12\varepsilon\mu_J+J$ of the leading form found in the proof of Proposition~\ref{prop:quad-wall}. Once $|\varepsilon|$ exceeds all these values, only the eigenvectors with $\mu_J=0$ contribute, that is, the zeros of $h$ on $\mathbb S^2$, and the degree is the sum of their indices. Generically, a quadratic $h$ has no zero-eigenvalue eigenvector. Adding $|v|^2c$ translates $h|_{\mathbb S^2}$ by $c\in\RR^3$, so Sard's theorem gives arbitrarily small shifts avoiding zero \cite[Chapter~II, \S1]{golubitskyguillemin}. If $V_h$ has only non-degenerate zeros, then sufficiently small $c$ preserve this property, since the change $\Pi_Jc$ tends to zero in $C^1$. Compactness makes the absence of zero-eigenvalue eigenvectors open, and the terminal degree is then zero. Theorem~\ref{thm:transition} shows a similar passage through a loss of properness for a radial perturbation. The escape to infinity at the threshold is particular to the degree $D=2$. For gradient perturbations of degree $D\geq3$, Subsection~\ref{subsec:anisotropy-reach} shows that the roots over an eigenvector of non-zero eigenvalue disappear in pairs, through a fold at finite radius.
\end{remark}

\begin{theorem}[Fourteen roots and a second degree transition]\label{thm:sphere-xyz}
For $q=x_0+x_1i+x_2j+x_3k$ and real $\varepsilon$, put
\[
G_\varepsilon(q)=q^2+1+\varepsilon\,(x_2x_3\,i+x_1x_3\,j+x_1x_2\,k).
\]
This map has real coordinate degree two. For $s=(s_1,s_2,s_3)\in\{-1,1\}^3$, write $\sigma=s_1s_2s_3$. The following holds:
\begin{enumerate}
\item For $0<|\varepsilon|<2\sqrt3$, the map is proper of Brouwer degree two and has exactly fourteen roots, all simple. The six roots $\pm i,\pm j,\pm k$ have local degree $-1$. The eight roots
\begin{equation}\label{eq:sphere-xyz-roots}
q=\frac{-\varepsilon\sigma+2(s_1i+s_2j+s_3k)}{\sqrt{12-\varepsilon^2}}
\end{equation}
have local degree $+1$.
\item For $|\varepsilon|>2\sqrt3$, the map is proper of Brouwer degree $-6$. Its roots are exactly $\pm i,\pm j,\pm k$, all simple and of local degree $-1$.
\item For $|\varepsilon|=2\sqrt3$, the map is not proper.
\end{enumerate}
\end{theorem}
\begin{proof}
The perturbation equals $\varepsilon h(\Ima q)$, where $h$ is the gradient of $f(v)=x_1x_2x_3$ for $v=x_1i+x_2j+x_3k$, so $D=2$. We first determine $\mathcal E(h)$. The eigenvector equations read $J_2J_3=\mu J_1$, $J_1J_3=\mu J_2$ and $J_1J_2=\mu J_3$. For $\mu=0$, at least two coordinates vanish, and $J$ is one of $\pm i,\pm j,\pm k$. For $\mu\neq0$, multiplication of the $a$-th equation by $J_a$ gives $\mu J_1^2=\mu J_2^2=\mu J_3^2=J_1J_2J_3$. The three squares then equal $\frac13$, so $J=(s_1i+s_2j+s_3k)/\sqrt3$ with $\mu=\sigma/\sqrt3$. Thus $\mathcal E(h)$ has fourteen points and $e(h)=7$.

The derivative $Dh(J)$ is the symmetric matrix with zero diagonal whose entry in position $(a,b)$ equals $J_c$, for $\{a,b,c\}=\{1,2,3\}$. At $J=i$ one has $\mu=0$, the tangent plane is spanned by $j$ and $k$, and $K_i$ has the matrix
\[
\begin{pmatrix}0&1\\1&0\end{pmatrix}
\]
of index $-1$. By Lemma~\ref{lem:eigen-antipodal} and the symmetry permuting the coordinates, the same index occurs at $\pm i,\pm j,\pm k$. At $J=(s_1i+s_2j+s_3k)/\sqrt3$, let $S=\operatorname{diag}(s_1,s_2,s_3)$ and let $U$ be the matrix whose entries all equal one. Then $Dh(J)=\frac{\sigma}{\sqrt3}S(U-\mathrm{Id})S$. A tangent vector has the form $\tau=Sw$ with $w$ orthogonal to $(1,1,1)$, and then $Dh(J)\tau=-\frac{\sigma}{\sqrt3}\tau$. Formula~\eqref{eq:eigen-derivative} gives $K_J=-\frac{2\sigma}{\sqrt3}\mathrm{Id}$, of index $+1$.

The eigenvalue $\mu=0$ satisfies $|\varepsilon\mu|<2$ for every $\varepsilon$, and~\eqref{eq:quad-wall-roots} gives the roots $\pm i,\pm j,\pm k$. For $|\mu|=1/\sqrt3$, the condition $|\varepsilon\mu|<2$ reads $|\varepsilon|<2\sqrt3$, and~\eqref{eq:quad-wall-roots} with $y=(1-\varepsilon^2/12)^{-1/2}$ gives~\eqref{eq:sphere-xyz-roots}. Proposition~\ref{prop:quad-wall} now proves all three assertions, the Brouwer degree being $8-6=2$ in the first range and $-6$ in the second.
\end{proof}

This theorem resembles Theorem~\ref{thm:transition} closely, and Remark~\ref{rem:quad-wall} describes the escape of the eight roots of positive local degree as $|\varepsilon|$ increases to $2\sqrt3$. The six roots at $\pm i,\pm j,\pm k$ remain in place for every $\varepsilon$, since $h$ vanishes on the coordinate axes. In the first range, the fourteen simple roots come close to the B\'ezout bound of sixteen for maps of degree two in Remark~\ref{rem:ambient}. The coercivity estimate in the proof of Proposition~\ref{prop:quad-wall} and Lemma~\ref{lem:stability} let this count persist on an open set of polynomial maps of degree at most two.

\subsection{Sharp bounds under a fixed perturbation degree}\label{subsec:spherical-fixed-degree}
The first reduced field of a polynomial perturbation of degree at most $D$ has the form
\[
V(J)=\Pi_J f(J),\qquad J\in\mathbb S^2
\]
for a polynomial map $f\colon\RR^3\to\RR^3$ of degree at most $D$. Indeed, restriction of $\Ima H$ to $x_0+y_0J$ does not increase its degree. Homogeneous imaginary perturbations provide examples attaining the largest possible number of non-degenerate zeros. A compactification of the sphere gives the same upper bound when several degrees occur together.

\begin{theorem}[A sharp bound for the first reduced field]\label{thm:sphere-bounded-count}
Let $D\geq1$ and let $f\colon\RR^3\to\RR^3$ have degree at most $D$. Assume that every zero of $V(J)=\Pi_Jf(J)$ on $\mathbb S^2$ is non-degenerate. Then
\begin{equation}\label{eq:sphere-bounded-count}
\#V^{-1}(0)\leq2(D^2+D+1).
\end{equation}
Equality occurs for a homogeneous map $f$ of degree $D$, which can be chosen as the gradient of a real homogeneous polynomial of degree $D+1$.

Consequently, for every normally non-degenerate spherical component, the largest local count obtained from a polynomial perturbation of degree at most $D$ whose first reduced field has only non-degenerate zeros is $2(D^2+D+1)$.
\end{theorem}

\begin{proof}
Consider the smooth complex projective quadric
\[
Q=\{[z_1:z_2:z_3:w]\in\mathbb P^3_{\CC}
\mid z_1^2+z_2^2+z_3^2=w^2\}.
\]
The coordinates $z_1,z_2,z_3$ never vanish simultaneously on $Q$, since $z_1=z_2=z_3=0$ would force $w=0$, leaving all four homogeneous coordinates zero. The resulting morphism has rank one everywhere and gives the exact sequence of vector bundles
\[
0\longrightarrow\mathcal O_Q(-1)
\xrightarrow{(z_1,z_2,z_3)}
\mathcal O_Q^{\oplus3}
\longrightarrow E
\longrightarrow0
\]
with $E$ of rank two. Homogenize the components of $f$ to degree $D$ by means of the coordinate $w$. On the chart $w=1$ they agree with $f$, and every real point of $Q$ lies in this chart, since $z_1^2+z_2^2+z_3^2=0$ has no non-zero real solution. Their classes modulo the line spanned by $z=(z_1,z_2,z_3)$ define a section $s_f$ of $E(D)=E\otimes\mathcal O_Q(D)$.

On the affine chart $w=1$, the vanishing of this section means that $f(z)$ is parallel to $z$. At a real point of this chart, $|z|=1$, so these equations are precisely $V(z)=0$. On $w=1$, the fiber of $E(D)$ is $\CC^3/\CC z$, identified with the complex tangent plane by $\operatorname{Id}-zz^{\mathsf T}$. At a real zero $J$ with $f(J)=\mu J$, the vertical derivative of $s_f$ is represented by the complexification of $(\Pi_JDf(J)-\mu\operatorname{Id})|_{T_J\mathbb S^2}$. This is the intrinsic derivative of $V$. Thus every non-degenerate real zero is a simple complex zero of $s_f$.

Write $\xi=c_1(\mathcal O_Q(1))$. The exact sequence gives
\[
c_1(E)=\xi,\qquad c_2(E)=\xi^2.
\]
For a bundle of rank two, tensoring by $\mathcal O_Q(D)$ gives
\[
c_2(E(D))
=c_2(E)+D\xi\,c_1(E)+D^2\xi^2
=(D^2+D+1)\xi^2.
\]
The integral $\int_Q\xi^2$ is the degree of the quadric, namely two, as one also sees from the isomorphism $Q\cong\mathbb P^1_{\CC}\times\mathbb P^1_{\CC}$, under which $\mathcal O_Q(1)$ corresponds to $\mathcal O(1,1)$. The top Chern number is therefore
\begin{equation}\label{eq:sphere-chern-count}
\int_Qc_2(E(D))=2(D^2+D+1).
\end{equation}

The bundle $E(D)$ is globally generated by triples of homogeneous polynomials of degree $D$. Indeed, at any $p\in Q$, choose a linear form $\ell$ with $\ell(p)\neq0$. The triples $c\ell^D$, with $c\in\CC^3$, span the fiber through the quotient map. A generic such section is transverse to the zero section. Since the bundle rank equals the complex dimension of $Q$, its zeros are isolated and simple, and compactness makes them finite. The top Chern class formula for a regular section gives their number as~\eqref{eq:sphere-chern-count} \cite[\S14.1]{fulton}. Choose disjoint neighborhoods of the non-degenerate real zeros of $V$, which are simple complex zeros of $s_f$. The implicit function theorem preserves one zero in each neighborhood under a sufficiently small generic perturbation. Each preserved point contributes one to the zero count of this transverse section, and every remaining zero has positive multiplicity. The number of non-degenerate real zeros of $V$ therefore cannot exceed~\eqref{eq:sphere-chern-count}. This perturbation argument applies regardless of positive-dimensional zero components of $s_f$, including those on the divisor $w=0$. This proves~\eqref{eq:sphere-bounded-count}.

For attainability, Cartwright and Sturmfels give $1+D+D^2=D^2+D+1$ projective eigenvector lines for a generic homogeneous map of degree $D$ in three variables \cite{cartwrightsturmfels}. For $D>1$, their tensor order is $D+1$, and their formula is $(D^3-1)/(D-1)$, while the sum gives three for $D=1$. Kozhasov constructs real homogeneous forms of degree $D+1$ whose gradients have all these eigenvector lines real and whose restrictions to the sphere are Morse \cite{kozhasov}. Each real line meets $\mathbb S^2$ in two points, and the Morse condition gives non-degenerate zeros of the projected gradient. This yields $2(D^2+D+1)$ zeros. Real eigenconfigurations of ternary forms are studied further in \cite{abo}. For $D=1$, a symmetric matrix with three distinct eigenvalues gives the same conclusion directly, Theorem~\ref{thm:sphere-xyz} gives an explicit example for $D=2$, and Proposition~\ref{prop:cubic-anisotropy} gives one for $D=3$.

Finally, for a chosen spherical component, take $H(q)=f(\Ima q)$ with such a homogeneous $f$. The first reduced field is $y_0^D\Pi_Jf(J)$, so it has the same zeros and indices. Theorem~\ref{thm:sphere-splitting} transfers the count to the nearby quaternionic roots.
\end{proof}

\begin{corollary}[Homogeneous imaginary perturbations]\label{cor:sphere-homogeneous-count}
Let $\Sigma$ be a spherical component with $\kappa\neq0$, and let $h\colon\Ima\HH\to\Ima\HH$ be homogeneous of degree $D\geq1$. Assume that every zero of $V_h$ is non-degenerate. For every sufficiently small $\varepsilon\neq0$, the perturbation $\varepsilon h(\Ima q)$ splits $\Sigma$ into exactly $2e(h)$ simple roots, where
\[
e(h)\leq D^2+D+1
\]
is odd. Equality occurs for gradients of real homogeneous forms of degree $D+1$. All these homogeneous local counts are congruent to two modulo four.
\end{corollary}

\begin{proof}
The field on the chosen sphere is $y_0^DV_h$. Theorem~\ref{thm:sphere-splitting} gives the correspondence with nearby roots. Lemma~\ref{lem:eigen-antipodal} gives the parity of $e(h)$, and Theorem~\ref{thm:sphere-bounded-count} gives the upper bound and its attainment.
\end{proof}

\begin{remark}[Signed counts and degenerate first fields]\label{rem:sphere-mixed}
The sharp bound~\eqref{eq:sphere-bounded-count} applies equally to homogeneous and mixed-degree first reduced fields. At the maximum, the counts of positive and negative quaternionic local degrees are
\[
(n_+,n_-)=
\begin{cases}
(D^2+D+2,\,D^2+D), & \kappa>0\\
(D^2+D,\,D^2+D+2), & \kappa<0.
\end{cases}
\]
These formulas follow from Theorem~\ref{thm:sphere-splitting}. For a first reduced field satisfying $V(-J)=V(J)$ or $V(-J)=-V(J)$, the zeros occur in pairs of equal index. The Poincar\'e--Hopf index sum of two forces an odd number of pairs, hence a count congruent to two modulo four. Homogeneity gives $h(-v)=(-1)^Dh(v)$ and ensures this parity. Mixed-degree fields can break it. For odd $\ell$, the perturbation $H_\ell$ of Proposition~\ref{prop:sphere-realization} contains a term of odd degree $\ell+2$ and a non-zero quadratic term, allowing counts divisible by four.

When the first reduced field has degenerate zeros, the correspondence in Theorem~\ref{thm:sphere-splitting} does not determine the complete nearby root count. Higher reduced terms, including~\eqref{eq:sphere-second-section}, can then influence the splitting. The bound proved here concerns non-degenerate first reduced fields and does not settle all such higher-order splittings.
\end{remark}

\subsection{Consequences for polynomials with one leading word}\label{subsec:spherical-one-word}
The theorem of Eilenberg and Niven \cite{en} guarantees a root for $q^d$ plus terms of lower degree. A one-sided polynomial of degree $d$ has at most $d$ isolated roots \cite{pogoruishapiro}. The spherical reduction shows that, among two-sided polynomials with the single leading word $q^d$, the number of simple roots can grow like $d^3$, and Subsection~\ref{subsec:en-stacked} raises this growth to $d^4$.

\begin{theorem}[Splitting the spheres of $q^d+1$]\label{thm:en-splitting}
Let $d=2m\geq2$ and $1\leq D\leq d-1$. Let $h$ be a homogeneous polynomial map of $\Ima\HH$ of degree $D$, with every unit eigenvector a non-degenerate zero of $V_h$. There is $\varepsilon_0>0$ such that, for $0<|\varepsilon|<\varepsilon_0$, the map
\[
G_\varepsilon(q)=q^d+1+\varepsilon h(\Ima q)
\]
has exactly $d\,e(h)$ roots, all simple. Exactly $m$ of them lie over each $J\in\mathcal E(h)$, and each of these has local degree $\operatorname{ind}_J(V_h)$. The map is proper of Brouwer degree $d$, and its root count persists on an open set of polynomial maps of degree at most $d$.
\end{theorem}
\begin{proof}
There is no real root, since $x^d+1>0$. For the auxiliary complex number $z=x+\mathrm iy$, one has $A(x,y)+\mathrm iB(x,y)=z^d+1$. The field $V_{\varepsilon h}=\varepsilon V_h$ has the same zeros and the same indices as $V_h$, the tangent plane being two-dimensional. By Lemma~\ref{lem:eigen-reduction}, applied to $\varepsilon h$, the roots over $J\in\mathcal E(h)$ with eigenvalue $\mu$ correspond to the zeros with $\Ima_{\mathrm i} z>0$ of
\[
\Phi_\tau(z)=z^d+1+\mathrm i\tau(\Ima_{\mathrm i} z)^D,\qquad \tau=\varepsilon\mu
\]
and the planar determinant in~\eqref{eq:eigen-det} is the Jacobian determinant of $\Phi_\tau$ as a real map.

Let $w_1,\dots,w_d$ be the roots of $w^d=-1$. None of them is real, exactly $m$ lie in the upper half-plane, and we put $\eta=\min_k|\Ima_{\mathrm i} w_k|>0$. The estimate $|z|^d\leq1+|\tau||z|^D$ confines the zeros of $\Phi_\tau$ with $|\tau|\leq1$ to a fixed disk, because $D<d$. As a real map, $z\mapsto z^d+1$ has Jacobian determinant $d^2|w_k|^{2d-2}=d^2$ at each $w_k$. The implicit function theorem and compactness give $\delta\in(0,\eta)$ and $\tau_0\in(0,1]$ with the following property. For $|\tau|<\tau_0$, the map $\Phi_\tau$ has exactly one zero in each disk $|z-w_k|<\delta$, with positive Jacobian determinant, and no zero outside these disks. The zeros near the $m$ roots in the upper half-plane have $\Ima_{\mathrm i} z>0$, and the others have $\Ima_{\mathrm i} z<0$.

The field $V_h$ has isolated zeros, so $h\neq0$. Every eigenvalue satisfies $|\mu|\leq\|h\|$, and the choice $\varepsilon_0=\tau_0/\|h\|$ makes the previous paragraph apply to every $J\in\mathcal E(h)$. Each unit eigenvector therefore carries exactly $m$ roots, and Lemma~\ref{lem:eigen-reduction} gives their simplicity and local degrees. The total is $m\cdot2e(h)=d\,e(h)$. The estimate $|G_\varepsilon(q)|\geq|q|^d-|\varepsilon|\,\|h\|\,|q|^D-1$ gives properness and the coercivity hypothesis of Lemma~\ref{lem:stability} with $N=d$. By the Poincar\'e--Hopf theorem, the indices over $\mathcal E(h)$ add up to two, so the signed count equals $2m=d$.
\end{proof}

With $D=d-1$ and $h$ as in Theorem~\ref{thm:sphere-bounded-count}, one has $e(h)=(d-1)^2+(d-1)+1=d^2-d+1$, and Theorem~\ref{thm:en-splitting} gives a polynomial of the form $q^d$ plus terms of degree less than $d$ with exactly $d(d^2-d+1)$ roots, all simple.

\begin{remark}[Odd exponents]\label{rem:en-odd}
The evenness of $d$ in Theorem~\ref{thm:en-splitting} is inherited from the radial family, and its proof uses it only to exclude a real root of $w^d=-1$. Let $d\geq3$ be odd, and let $D$ and $h$ be as in Theorem~\ref{thm:en-splitting}. For every sufficiently small $\varepsilon\neq0$, the map $q^d+1+\varepsilon h(\Ima q)$ then has exactly $(d-1)e(h)+1$ roots, all simple, and Brouwer degree $d$. Indeed, $(d-1)/2$ roots of $w^d=-1$ lie in the upper half-plane, and the argument of Theorem~\ref{thm:en-splitting} places $(d-1)/2$ roots over each of the $2e(h)$ unit eigenvectors. The planar map $\Phi_\tau$ agrees with $z^d+1$ on the real axis, so its unique zero near $-1$ is $-1$ itself, and no root with $y>0$ arises there. The only real root is $q=-1$. At this point the derivative is block diagonal with respect to $\RR\oplus\Ima\HH$, with blocks $d$ and $d\,\mathrm{Id}+\varepsilon Dh(0)$, so the root is simple for small $\varepsilon$ and has local degree $+1$. For $D=d-1$ the count is $(d-1)(d^2-d+1)+1$ for odd $d\geq3$. Together with Theorem~\ref{thm:en-splitting}, this gives cubic growth for every degree $d\geq2$. The stacking formula~\eqref{eq:en-stacked-g} uses the integer exponent $m=d/2$. For odd $d$, the radial equations $\pm u^{d/2}+g(u)=0$ lie outside its polynomial factorization argument.
\end{remark}

\subsection{Several spheres at once}\label{subsec:en-stacked}
The spheres of $q^d+1$ lie on a single circle of the complex slice. A lower-degree radial term can create many more spheres, each of which splits as in Theorem~\ref{thm:en-splitting}. This raises the growth of the root count from $d^3$ to $d^4$.

\begin{theorem}[Several spheres split at once]\label{thm:en-stacked}
Let $d=2m\geq2$ and $1\leq D\leq d-1$. Let $h$ be a homogeneous polynomial map of $\Ima\HH$ of degree $D$, with every unit eigenvector a non-degenerate zero of $V_h$. Choose numbers $0<r_1<\cdots<r_m$ and put
\begin{equation}\label{eq:en-stacked-g}
g(u)=u^m-\prod_{j=1}^m(u-r_j),\qquad G_\varepsilon(q)=q^d+g(|q|^2)+\varepsilon h(\Ima q).
\end{equation}
The polynomial $g$ has degree $m-1$, so $g(|q|^2)$ has real coordinate degree $d-2$, while $\varepsilon h(\Ima q)$ has degree $D\leq d-1$. Thus $G_\varepsilon$ is $q^d$ plus terms of degree less than $d$, and it belongs to the class of Eilenberg and Niven. There is $\varepsilon_0>0$ such that, for $0<|\varepsilon|<\varepsilon_0$, the following holds:
\begin{enumerate}
\item The map $G_\varepsilon$ is proper of Brouwer degree $d$, and it has at least $2m^2e(h)$ simple roots, with at least $m^2$ over each unit eigenvector of $h$.
\item Suppose that the positive roots of $2u^m-\prod_{j=1}^m(u-r_j)$ are simple, and let $k_+$ be their number. Then $G_\varepsilon$ has exactly $2e(h)\bigl(m^2+(m-1)k_+\bigr)+2k_+$ roots, all simple. For each such $\varepsilon$, this exact count persists on an open neighborhood of $G_\varepsilon$ in the coefficient space of polynomial maps of degree at most $d$.
\end{enumerate}
\end{theorem}
\begin{proof}
The coefficient of $u^{m-1}$ in $g$ is $r_1+\cdots+r_m>0$, so $g$ has degree $m-1$ and $g(|q|^2)$ has real coordinate degree $d-2$. For $z\in\CC$ and real $\tau$ put
\begin{equation}\label{eq:en-stacked-planar}
\Phi_\tau(z)=z^d+g(|z|^2)+\mathrm i\tau(\Ima_{\mathrm i} z)^D.
\end{equation}
The map $F_0(q)=q^d+g(|q|^2)$ satisfies the hypothesis of Lemma~\ref{lem:eigen-reduction}, with $A(x,y)+\mathrm iB(x,y)=\Phi_0(x+\mathrm iy)$. We apply the lemma to $\varepsilon h$, whose unit eigenvectors are those of $h$ and whose field $\varepsilon V_h$ has the same indices. The roots of $G_\varepsilon$ over $J\in\mathcal E(h)$ with eigenvalue $\mu$ then correspond to the zeros of $\Phi_{\varepsilon\mu}$ in the upper half-plane, and the planar determinant in~\eqref{eq:eigen-det} is the Jacobian determinant of $\Phi_{\varepsilon\mu}$ as a real map. The real roots of $G_\varepsilon$ are the real zeros of $\Phi_0$.

We first locate the zeros of $\Phi_0$. The value $g(0)=(-1)^{m+1}r_1\cdots r_m$ is non-zero, so $z=0$ is not a zero. For $z=\sqrt u\,e^{\mathrm i\theta}$ with $u>0$ one has $\Phi_0(z)=u^me^{\mathrm id\theta}+g(u)$. A zero therefore has $e^{\mathrm id\theta}=\sigma\in\{-1,1\}$ and $\sigma u^m+g(u)=0$. For $\sigma=-1$ this equation reads $\prod_j(u-r_j)=0$, and the zeros in the upper half-plane are the $m^2$ points
\begin{equation}\label{eq:en-stacked-zeros}
z_{j,l}=\sqrt{r_j}\,e^{\mathrm i(2l+1)\pi/d},\qquad 1\leq j\leq m,\qquad 0\leq l\leq m-1.
\end{equation}
No zero with $\sigma=-1$ is real, since $x^d\geq0$ for real $x$. For $\sigma=1$ the equation reads $2u^m-\prod_j(u-r_j)=0$. Each of its positive roots $u$ gives the $m-1$ zeros $\sqrt u\,e^{2\pi\mathrm il/d}$ with $1\leq l\leq m-1$ in the upper half-plane, together with the two real zeros $\pm\sqrt u$. The simplicity hypothesis holds for an open dense set of ordered radii, and in particular for distinct radii sufficiently close to a common $r>0$. Indeed, $2u^m-(u-r)^m$ and its derivative have no common zero, so the discriminant is a non-zero polynomial in the radii.

A zero $z_*\neq0$ of $\Phi_0$, with $u_*=|z_*|^2$ and sign $\sigma$, is non-degenerate exactly when $u_*$ is a simple root of $\sigma u^m+g(u)$. Indeed, write $\Phi_0=T\circ Q$ with $Q(z)=z^d$ and $T(w)=w+\varphi(|w|)$, where $\varphi(t)=g(t^{2/d})$ is smooth for $t>0$. The Jacobian determinant of $Q$ at $z_*$ equals $d^2|z_*|^{2d-2}>0$. At the real point $w_*=\sigma u_*^m$, the derivative of $T$ in the real and imaginary coordinates is $\operatorname{diag}(1+\sigma\varphi'(t_*),1)$ with $t_*=u_*^m$, as in the proof of Lemma~\ref{lem:sphere-normal}. Moreover $1+\sigma\varphi'(t)=\sigma\frac{d}{dt}\bigl(\sigma t+\varphi(t)\bigr)$, and the substitution $t=u^m$, which has positive derivative for $u>0$, turns $\sigma t+\varphi(t)$ into $\sigma u^m+g(u)$. The roots $r_j$ are simple, so the $m^2$ zeros~\eqref{eq:en-stacked-zeros} are non-degenerate.

There is a constant $C$ depending only on $g$ such that
\[
|\Phi_\tau(z)|\geq|z|^d-C\bigl(1+|z|^{d-2}\bigr)-|\tau|\,|z|^D.
\]
Since $D<d$, the zeros of $\Phi_\tau$ with $|\tau|\leq1$ lie in a fixed disk. The implicit function theorem gives $\tau_0\in(0,1]$ and, for each point~\eqref{eq:en-stacked-zeros}, a neighborhood in the upper half-plane in which $\Phi_\tau$ has exactly one zero for $|\tau|<\tau_0$, with Jacobian determinant of the same sign as for $\Phi_0$. The field $V_h$ has isolated zeros, so $h\neq0$, and every eigenvalue satisfies $|\mu|\leq\|h\|$. Put $\varepsilon_0=\tau_0/\|h\|$. For $0<|\varepsilon|<\varepsilon_0$, each unit eigenvector therefore carries at least $m^2$ roots, and each of them is simple by Lemma~\ref{lem:eigen-reduction}, both determinants in~\eqref{eq:eigen-det} being non-zero. Distinct eigenvectors give distinct roots, and $\mathcal E(h)$ has $2e(h)$ points by Lemma~\ref{lem:eigen-antipodal}. The homotopy $q^d+s\bigl(g(|q|^2)+\varepsilon h(\Ima q)\bigr)$, for $0\leq s\leq1$, is proper by an estimate of the same form, so the Brouwer degree equals $d$. This proves the first assertion.

Under the hypothesis of the second assertion, every zero of $\Phi_0$ is non-degenerate, and they are finite in number since they lie in a fixed disk. After shrinking $\tau_0$, compactness shows that $\Phi_\tau$ has no zero outside small neighborhoods of these points, each neighborhood containing exactly one zero. The neighborhoods of non-real zeros may be taken in the open half-plane containing them. On the real axis $\Phi_\tau=\Phi_0$, since $\Ima_{\mathrm i} z=0$ there, so each real zero of $\Phi_0$ is the unique zero of $\Phi_\tau$ in its neighborhood. The zeros of $\Phi_\tau$ in the upper half-plane are therefore exactly the $m^2+(m-1)k_+$ perturbed zeros, and each unit eigenvector carries exactly this many roots. The real roots of $G_\varepsilon$ are the $2k_+$ points $x_*=\pm\sqrt u$. At such a point, the derivative of $G_\varepsilon$ is block diagonal with respect to $\RR\oplus\Ima\HH$. Its real block is $\partial_xA(x_*,0)$, and its imaginary block is $d\,x_*^{d-1}\mathrm{Id}+\varepsilon Dh(0)$, since the derivative of $g(|q|^2)$ in an imaginary direction vanishes at a real point. The symmetry $\Phi_0(\overline z)=\overline{\Phi_0(z)}$ makes the planar Jacobian at $x_*$ diagonal, with $\partial_yB(x_*,0)=d x_*^{d-1}\neq0$, so planar non-degeneracy gives $\partial_xA(x_*,0)\neq0$. The imaginary block is invertible for small $\varepsilon$, since $x_*\neq0$. Adding the two contributions gives the stated count. For fixed $\varepsilon$, the estimate $|G_\varepsilon(q)|\geq|q|^d-O(|q|^{d-1}+1)$ gives the coercivity hypothesis of Lemma~\ref{lem:stability} with $N=d$, which proves persistence.
\end{proof}

\begin{remark}[Local degrees in Theorem~\ref{thm:en-stacked}]\label{rem:en-stacked-signs}
Under the hypothesis of the second assertion, the local degrees can be read from the proof. Put $p(u)=2u^m-\prod_{j=1}^m(u-r_j)$. Over $J\in\mathcal E(h)$, the root near a zero of $\Phi_0$ with $|z|^2=u_*$ and sign $\sigma$ has local degree $s\,\operatorname{ind}_J(V_h)$, where $s$ is the sign of the planar determinant, that is, the sign of $\sigma\frac{d}{du}\bigl(\sigma u^m+g(u)\bigr)$ at $u_*$. For $\sigma=-1$ and $u_*=r_j$ one finds $s=(-1)^{m-j}$. For $\sigma=1$ one finds $s=\sgn p'(u_*)$, and these signs alternate along the positive roots of $p$, ending with $+1$ at the largest. At a real root $x_*=\pm\sqrt{u_*}$, one has $A(x,0)=p(x^2)$, so the real block is $2x_*p'(u_*)$, while the imaginary block has the sign of $x_*^{3(d-1)}$, that is, of $x_*$, for small $\varepsilon$. The two real roots over $u_*$ have therefore the common local degree $\sgn p'(u_*)$.

Since $p(0)=(-1)^{m+1}\prod_jr_j$ and $p(u)\to+\infty$, $k_+$ is even for odd $m$ and odd for even $m$, so the sum of $\sgn p'(u_*)$ over the positive roots is $0$ and $1$, respectively. The sum of the signs $s$ over the zeros in the upper half-plane equals $m$ for $m$ odd and $m-1$ for $m$ even, and the real roots contribute $0$ and $2$ respectively. Since the indices over $\mathcal E(h)$ add up to two, the signed total is $2m=d$ in both cases, in agreement with the Brouwer degree.

For the polynomial~\eqref{eq:en-quartic-example} of Corollary~\ref{cor:en-quartic} with $h=\nabla W_c$ and $W_c$ given by~\eqref{eq:cubic-energy}, each unit eigenvector carries three roots with $s=+1$ and two with $s=-1$. By Proposition~\ref{prop:cubic-anisotropy} there are fourteen eigenvectors of index $+1$ and twelve of index $-1$. Together with the two real roots, this gives $14\cdot3+12\cdot2+2=68$ roots of local degree $+1$ and $14\cdot2+12\cdot3=64$ roots of local degree $-1$, with difference four.
\end{remark}

\begin{corollary}[Quartic growth in the one-word class]\label{cor:en-quartic}
For every even $d\geq2$, there is a two-sided polynomial of the form $q^d$ plus terms of degree less than $d$ with exactly $d(d-1)(d^2-2d+3)$ roots, all simple, and every polynomial map of degree at most $d$ sufficiently close to it has the same number of roots, all simple. For $d=4$, take $h=\nabla W_c$ with $W_c(x)=x_1^2x_2^2+x_2^2x_3^2+x_3^2x_1^2$, whose eigenvectors are described in Proposition~\ref{prop:cubic-anisotropy}. For every $0<|\varepsilon|\leq3/40$, the polynomial
\begin{equation}\label{eq:en-quartic-example}
q^4+3|q|^2-2+\varepsilon h(\Ima q)
\end{equation}
has exactly $132$ roots, all simple. This count persists on an open set of polynomial maps of degree at most four.
\end{corollary}
\begin{proof}
Take $D=d-1$ and $h$ as in Theorem~\ref{thm:sphere-bounded-count}, with $e(h)=d^2-d+1$. Choose $r_j=4^j$ in Theorem~\ref{thm:en-stacked} and put $p(u)=2u^m-\prod_j(u-r_j)$. For $1\leq k<m$ with $m-k$ even, the factors $1-r_j/(2r_k)$ are positive for $j\leq k$ and negative for $j>k$. The latter occur in an even number, so their product is positive. The first group has product at least $\frac12(1-\sum_{s\geq1}(2\cdot4^s)^{-1})=5/12$. The factors with $j=k+1,k+2$ have product $(-1)(-7)=7$, and every remaining factor has absolute value greater than one. Hence $\prod_j(1-r_j/(2r_k))\geq35/12>2$, giving $p(2r_k)<0$. Since $r_k<2r_k<r_{k+1}$ and $p(r_k),p(r_{k+1})>0$, there is a sign change in each of $(r_k,2r_k)$ and $(2r_k,r_{k+1})$.

For even $m$, $p(0)<0<p(r_1)$ gives one further root in $(0,r_1)$. These intervals yield $m-1$ distinct positive roots. Moreover, $p(0)=(-1)^{m+1}\prod_jr_j$ has the opposite sign from $p(u)$ as $u\to-\infty$, giving a negative root. These $m$ distinct roots exhaust the degree and are simple. The second assertion of Theorem~\ref{thm:en-stacked} now gives $2e(h)\bigl(m^2+(m-1)^2\bigr)+2(m-1)=d(d-1)(d^2-2d+3)$ simple roots. The persistence assertion is part of that theorem.

For $d=4$, take $m=2$, $r_1=1$ and $r_2=2$. Then $g(u)=u^2-(u-1)(u-2)=3u-2$, and $2u^2-(u-1)(u-2)=u^2+3u-2$ has the single positive root $(\sqrt{17}-3)/2$, which is simple. Proposition~\ref{prop:cubic-anisotropy} in Section~\ref{sec:anisotropy} computes the eigenvalues $0,1,4/3$ and the counts $(14,12)$ of positive and negative tangent indices for $h=\nabla W_c$, giving $e(h)=13$. The second assertion of Theorem~\ref{thm:en-stacked} with $e(h)=13$ and $k_+=1$ gives $2\cdot13\cdot5+2=132$ simple roots. Since every eigenvalue satisfies $|\mu|\leq4/3$, the restriction $|\varepsilon|\leq3/40$ gives $|\tau|=|\varepsilon\mu|\leq1/10$. For $\tau>0$, the planar imaginary equation gives $y^2=4x^3/(4x-\tau)$. Substitution into the real equation and multiplication by $(4x-\tau)^2$ give the degree-six polynomial
\begin{equation}\label{eq:en-quartic-planar-polynomial}
\begin{aligned}
P_\tau(x)&=-32x^2(2x^2-1)(x^2-1)\\
&\quad+4\tau x(4x^4-9x^2+4)+\tau^2(x^4+3x^2-2).
\end{aligned}
\end{equation}
For $0<\tau\leq1/10$, the terms independent of $\tau$ at $x=\pm1/2,\pm4/5,\pm6/5$ are respectively $-3$, $32256/15625>2$ and $-595584/15625<-38$. At these six points, the coefficients of $\tau$ and $\tau^2$ have absolute values at most $4$ and $5$, so their combined contribution is at most $4\tau+5\tau^2\leq9/20<1/2$ and preserves those signs. The remaining evaluations are
\begin{equation}\label{eq:en-quartic-sign-tests}
\begin{aligned}
P_\tau(0)&=-2\tau^2,& P_\tau(\tau/4)&=\tau^6/256,\\
P_\tau(\tau/2)&=-\tau^2\bigl(2-\tfrac94\tau^2+\tfrac7{16}\tau^4\bigr).
\end{aligned}
\end{equation}
The last bracket is at least $791/400>0$. These signs give six distinct real roots, exhausting degree six, so every root of $P_\tau$ is simple. Exactly five give $y^2>0$, the root in $(0,\tau/4)$ being excluded. At a retained planar solution, the derivative of the imaginary equation with respect to $y$ equals $2(\tau-4x)y^2=-8x^3\neq0$. Differentiating the eliminated real equation gives
\begin{equation}\label{eq:en-quartic-planar-jacobian}
\det D\Phi_\tau(x+\mathrm i y)=-\frac{8x^3P_\tau'(x)}{(4x-\tau)^2}\neq0.
\end{equation}
Thus every retained planar root is simple. Reflection in $x$ treats $\tau<0$, and the five upper half-plane zeros of $\Phi_0$ treat $\tau=0$. Thus each of the $26$ eigenvectors carries five roots, together with the two real roots. Persistence follows from the coercivity estimate in the proof of Theorem~\ref{thm:en-stacked} and Lemma~\ref{lem:stability}.
\end{proof}

\subsection{The quadratic one-word class}\label{subsec:en-quadratic}
For $d=2$, Theorems~\ref{thm:en-splitting} and~\ref{thm:en-stacked} both give six roots, and Proposition~\ref{prop:sphere-linear} realizes them for every $\varepsilon\neq0$. The quadratic class admits a complete answer, which shows that the spherical constructions do not reach the optimum there.

\begin{proposition}[The quadratic one-word class]\label{prop:en-quadratic}
For every real affine map $L\colon\HH\to\HH$, each fiber of $q\mapsto q^2+L(q)$ has at most eight simple roots, since a target is absorbed into the constant term of $L$. The bound is attained in the zero fiber by
\begin{equation}\label{eq:en-quadratic-example}
G_*(q)=q^2-(2x_0+x_1+x_2+2x_3)-2(x_0+x_3)\,i+2(x_1-x_3)\,j-3x_2\,k
\end{equation}
which has eight simple roots, five of local degree $+1$ and three of local degree $-1$. Among them are $0$, $1-\frac12i-j-\frac32k$ and $\pm\sqrt{3/2}-j\mp\sqrt{3/2}\,k$.
\end{proposition}
\begin{proof}
Write $q=x+v$ and decompose
\[
L(x+v)=\bigl(\alpha x+\langle\beta,v\rangle+c_0\bigr)+\bigl(xa+Mv+c\bigr)
\]
with real $\alpha,c_0$, vectors $a,\beta,c\in\Ima\HH$ and a linear map $M$ of $\Ima\HH$. The imaginary part of the equation reads $(2x+M)v=-(xa+c)$. Put $\delta(x)=\det(2x+M)$ and $W(x)=-\operatorname{adj}(2x+M)(xa+c)$. At a root with $\delta(x)\neq0$, one has $v=W(x)/\delta(x)$, and the real part of the equation becomes $p(x)=0$ with
\[
p=(x^2+\alpha x+c_0)\,\delta^2-|W|^2+\delta\,\langle\beta,W\rangle.
\]
The polynomial $\delta$ has degree three with leading coefficient eight, and $W$ has degree at most three. Hence $p$ has degree eight with leading coefficient $64$ for every affine map $L$. In particular $p$ never vanishes identically, and there are at most eight roots with $\delta(x)\neq0$, since distinct roots of this kind have distinct real parts. At $\delta(x)=0$, the condition $xa+c\in\operatorname{im}(2x+M)$ makes the imaginary equation solvable, with an affine space of candidates. The real equation leaves an empty set, a point or a sphere in that space, and positive-dimensional spheres give non-isolated roots. Matrices $M$ with three distinct complex eigenvalues form a dense set. For such an $M$, the real zeros of $\delta$ are the points $x=-\mu/2$ for the at most three real eigenvalues $\mu$ of $M$, each of them simple. For each fixed pair $(M,a)$ and each of these zeros, the solvability condition confines $c$ to a proper affine plane, since $2x+M$ has rank two. The affine maps $L$ for which $M$ has three distinct complex eigenvalues and $c$ avoids these planes therefore form a dense set, and each of them has at most eight roots. Simple roots are stable under small perturbations, whether or not $\delta$ vanishes at them. If some $L$ had nine simple roots, then the implicit function theorem would preserve them for every sufficiently close affine map, and in particular for one in the dense set, contrary to that bound. Thus every $L$ has at most eight simple roots.

Write $G(q)=q^2+L(q)$ and $g_0(x,v)=x^2-|v|^2+\alpha x+\langle\beta,v\rangle+c_0$ for its real component. For $\delta(x)\neq0$, substitution gives $g_0(x,W(x)/\delta(x))=p(x)/\delta(x)^2$. At a root, the Schur complement of $2x+M$ therefore gives
\begin{equation}\label{eq:en-quadratic-schur}
\det DG(q)=\delta(x)\left.\frac{d}{dx}\left(\frac{p(x)}{\delta(x)^2}\right)\right|_{p(x)=0}=\frac{p'(x)}{\delta(x)}.
\end{equation}
For~\eqref{eq:en-quadratic-example}, one has $\alpha=-2$, $c_0=0$, $c=0$, $\beta=-i-j-2k$, $a=-2i$ and
\[
M=\begin{pmatrix}0&0&-2\\2&0&-2\\0&-3&0\end{pmatrix},\qquad \delta(x)=8x^3-12x+12.
\]
A direct computation gives
\[
p(x)=16\,x\,(x-1)(2x^2-3)(2x^4-2x^3-9x^2+5x+3).
\]
Put
\[
f(x)=2x^4-2x^3-9x^2+5x+3.
\]
Its values satisfy
\[
f(-2)=5,\qquad f(-7/4)=-491/128,\qquad f(-1)=-7,\qquad f(0)=3
\]
and
\[
f(1)=-1,\qquad f(2)=-7,\qquad f(3)=45.
\]
Hence $f$ has a root in each of the intervals
\[
(-2,-7/4),\qquad (-1,0),\qquad (0,1),\qquad (2,3).
\]
These four roots are distinct and exhaust its degree, so all are simple. They are distinct from $0$, $1$ and $\pm\sqrt{3/2}$. Thus $p$ has eight distinct real roots.

The polynomial $\delta(x)=8x^3-12x+12$ has its two critical values equal to $12\pm4\sqrt2$, both positive. It therefore has a unique real zero. Since $\delta(-7/4)=-79/8$ and $\delta(-3/2)=3$, that zero lies in $(-7/4,-3/2)$, which contains no root of $p$. Consequently every root of $p$ gives exactly one quaternionic root through $v=W(x)/\delta(x)$. For completeness, the elimination vector is
\[
W(x)=\bigl(8x^3-12x,\,-8x^2,\,-12x\bigr).
\]

The leading coefficient of $p$ is positive, so the signs of $p'$ at its eight ordered roots alternate, beginning with $-1$. The value of $\delta$ is negative at the first root and positive at the other seven. Thus the signs of $p'/\delta$ are
\[
+1,\ +1,\ -1,\ +1,\ -1,\ +1,\ -1,\ +1.
\]
All eight quaternionic roots are simple, with five positive local degrees and three negative local degrees. Finally, $G_*(q)=q^2+O(|q|)$ is properly homotopic to $q^2$, so its Brouwer degree is two. The estimate $|G_*(q)|\geq|q|^2-O(|q|+1)$ and Lemma~\ref{lem:stability} give an open set of maps of degree at most two with exactly eight simple roots.
\end{proof}

The stacked construction retains the separation into eigenvector and planar equations of Lemma~\ref{lem:eigen-reduction}. General lower-degree terms can couple the real and imaginary coordinates, and the affine example~\eqref{eq:en-quadratic-example} exceeds the spherical count for $d=2$. For even $d$, let $M_d$ be the largest simple root count in the one-word class. Corollary~\ref{cor:en-quartic} and Remark~\ref{rem:ambient} give $d(d-1)(d^2-2d+3)\leq M_d\leq d^4$, so $M_d/d^4\to1$. The exact maximum for even $d\geq4$ remains open here.

\subsection{Scope of the spherical reduction}\label{subsec:spherical-scope}
The normal non-degeneracy hypothesis is essential to the reduction. When $\kappa=0$, the matrix $L_\Sigma$ is singular and the two normal variables cannot be eliminated by the implicit function theorem. Such components require additional parameters or a different local model. The unbounded spherical fibers in~\eqref{eq:threshold-zero-fiber} provide explicit examples at the threshold of Theorem~\ref{thm:transition}.

For $\kappa\neq0$, Proposition~\ref{prop:sphere-realization} realizes every positive even simple local count when the perturbation degree is unrestricted. Among perturbations of degree at most $D$ whose first reduced field has only non-degenerate zeros, Theorem~\ref{thm:sphere-bounded-count} gives the sharp maximum $2(D^2+D+1)$. Homogeneous imaginary perturbations attain this maximum and impose the additional antipodal restriction that their local counts are congruent to two modulo four. Mixed-degree perturbations obey the same upper bound without that parity restriction, as explained in Remark~\ref{rem:sphere-mixed}.

\section{An anisotropy dictionary}\label{sec:anisotropy}
The eigenvector reduction of Section~\ref{sec:eigen} has a familiar counterpart in physics. The direction of a classical spin, of a magnetization or of a molecular axis is a point $J\in\mathbb S^2$, and the anisotropy of the surrounding medium is often modeled by an energy that is the restriction to $\mathbb S^2$ of a real homogeneous polynomial. Its equilibrium directions are the critical points of this energy. They are the unit eigenvectors of the gradient perturbation, and their Morse types determine the tangent indices in the local-degree formulas~\eqref{eq:sphere-local-sign} and~\eqref{eq:eigen-det}. Three classical energies realize the equality case of Theorem~\ref{thm:sphere-bounded-count} for $D=1,2,3$. The quaternionic equations remain mathematical objects, and the correspondence transfers equilibrium counts, tangent indices and energy thresholds. We retain the basis $e_1=i$, $e_2=j$, $e_3=k$ of Subsection~\ref{subsec:spherical-realization}.

\subsection{Anisotropy energies and gradient perturbations}\label{subsec:anisotropy-energies}
Let $W$ be a real homogeneous polynomial of degree $D+1\geq2$ on $\Ima\HH=\RR^3$, and put
\begin{equation}\label{eq:anisotropy-energy}
E=W|_{\mathbb S^2},\qquad h=\nabla W.
\end{equation}
We call $E$ the anisotropy energy and $h$ the associated gradient perturbation. A critical point of $E$ is an equilibrium direction. It is a local minimum, a local maximum or a saddle when it is non-degenerate, since $\mathbb S^2$ has dimension two. The energy $E$ is a Morse function when all its critical points are non-degenerate.

\begin{lemma}[The anisotropy dictionary]\label{lem:anisotropy}
Let $W$, $E$ and $h$ be as in~\eqref{eq:anisotropy-energy}.
\begin{enumerate}
\item The field $V_h$ is the spherical gradient of $E$. Its zeros $\mathcal E(h)$ are the equilibrium directions, and the eigenvalue of $J\in\mathcal E(h)$ equals $\mu_J=(D+1)E(J)$.
\item At $J\in\mathcal E(h)$, the endomorphism $K_J$ of~\eqref{eq:eigen-derivative} is the Hessian of $E$, viewed as an endomorphism of $T_J\mathbb S^2$ through the inner product. Thus $J$ is a non-degenerate zero of $V_h$ exactly when it is a non-degenerate critical point of $E$, and then $\operatorname{ind}_J(V_h)$ equals $+1$ at a local minimum or maximum and $-1$ at a saddle.
\item Suppose that $E$ is a Morse function with $n_{\min}$ minima, $n_{\max}$ maxima and $n_{\mathrm{sad}}$ saddles. Then
\begin{equation}\label{eq:anisotropy-morse}
n_{\min}+n_{\max}-n_{\mathrm{sad}}=2,\qquad n_{\min}+n_{\max}+n_{\mathrm{sad}}\leq2(D^2+D+1).
\end{equation}
\end{enumerate}
\end{lemma}
\begin{proof}
The spherical gradient of the restriction of $W$ is the tangential component of $\nabla W$, which is $V_h$. The Euler identity gives $\langle\nabla W(J),J\rangle=(D+1)W(J)$, and at a zero of $V_h$ the vector $h(J)$ is parallel to $J$, so its eigenvalue is $(D+1)E(J)$. For the second assertion, let $J\in\mathcal E(h)$ and $\tau\in T_J\mathbb S^2$ be a unit vector. The great circle $\gamma(s)=\cos s\,J+\sin s\,\tau$ is a geodesic with $\gamma'(0)=\tau$ and $\gamma''(0)=-J$, so
\[
\frac{d^2}{ds^2}\Big|_{s=0}W(\gamma(s))=\langle D^2W(J)\tau,\tau\rangle-\langle\nabla W(J),J\rangle=\langle(Dh(J)-\mu_J)\tau,\tau\rangle.
\]
The Hessian of $E$ at the critical point $J$ is the symmetric form determined by these second derivatives, and its matrix in an orthonormal basis of $T_J\mathbb S^2$ is that of $K_J=(\Pi_JDh(J)-\mu_J)|_{T_J\mathbb S^2}$, since $Dh=D^2W$ is symmetric. The index of a non-degenerate zero of a gradient field is the sign of the Hessian determinant, which is positive for a definite Hessian and negative for an indefinite one in dimension two. The third assertion follows from Lemma~\ref{lem:eigen-antipodal}, in the form of the Poincar\'e--Hopf theorem on $\mathbb S^2$, and from Theorem~\ref{thm:sphere-bounded-count}.
\end{proof}

Table~\ref{tab:anisotropy} summarizes the correspondence. The bound in~\eqref{eq:anisotropy-morse} is the number of equilibrium directions permitted to a Morse energy of degree $D+1$, and it is of course a restatement of the eigenvector count of Cartwright and Sturmfels \cite{cartwrightsturmfels}. The quaternionic side adds the local degrees of the roots, and the thresholds of Proposition~\ref{prop:quad-wall}.

\begin{table}[ht]
\caption{The correspondence of Lemma~\ref{lem:anisotropy} and Proposition~\ref{prop:anisotropy-roots}, with the thresholds for $D\geq3$ given in~\eqref{eq:anisotropy-higher-threshold}.}
\label{tab:anisotropy}
\centering
\begin{tabular}{@{}ll@{}}
\toprule
Quaternionic side & Anisotropy side\\
\midrule
unit eigenvector $J\in\mathcal E(h)$ & equilibrium direction of $E$\\
eigenvalue $\mu_J$ & $(D+1)E(J)$\\
endomorphism $K_J$ & Hessian of $E$ at $J$, as an endomorphism\\
index $\operatorname{ind}_J(V_h)=\pm1$ & extremum or saddle\\
Poincar\'e--Hopf sum & $n_{\min}+n_{\max}-n_{\mathrm{sad}}=2$\\
bound $2(D^2+D+1)$ & largest number of Morse equilibria\\
threshold $|\varepsilon\mu_J|=2$ for $D=2$ & energy level $|E(J)|=2/(3|\varepsilon|)$\\
\bottomrule
\end{tabular}
\end{table}

\subsection{Roots over equilibrium directions}\label{subsec:anisotropy-roots}
The results of Section~\ref{sec:eigen} now read as statements about equilibria.

\begin{proposition}[Roots and equilibrium directions]\label{prop:anisotropy-roots}
Let $W$, $E$ and $h$ be as in~\eqref{eq:anisotropy-energy}, with $E$ a Morse function.
\begin{enumerate}
\item Let $d=2m$ with $D\leq d-1$. For every sufficiently small $\varepsilon\neq0$, the map $q^d+1+\varepsilon\nabla W(\Ima q)$ has exactly $m$ roots over each equilibrium direction, and no other roots. Those over minima and maxima have local degree $+1$, and those over saddles have local degree $-1$. It has therefore $m(n_{\min}+n_{\max})$ roots of local degree $+1$ and $m\,n_{\mathrm{sad}}$ roots of local degree $-1$.
\item Let $D=2$ and $\varepsilon\neq0$. The roots of $q^2+1+\varepsilon\nabla W(\Ima q)$ lie over the equilibrium directions $J$ with $|E(J)|<2/(3|\varepsilon|)$, one over each, with local degree $+1$ over an extremum and $-1$ over a saddle. The map is proper exactly when no equilibrium direction has $|E(J)|=2/(3|\varepsilon|)$, and its Brouwer degree is then the number of extrema minus the number of saddles among the equilibria with $|E(J)|<2/(3|\varepsilon|)$.
\end{enumerate}
\end{proposition}
\begin{proof}
By Lemma~\ref{lem:anisotropy}, every unit eigenvector of $h$ is a non-degenerate zero of $V_h$, with index given by its Morse type. The first assertion is then Theorem~\ref{thm:en-splitting}. For $D=2$ the eigenvalue is $\mu_J=3E(J)$, so the condition $|\varepsilon\mu_J|<2$ of Proposition~\ref{prop:quad-wall} reads $|E(J)|<2/(3|\varepsilon|)$, and the second assertion is that proposition.
\end{proof}

The second assertion has a simple reading. As the coupling $|\varepsilon|$ grows, the roots over the equilibria of largest energy in absolute value are the first to escape to infinity. Once $|\varepsilon|$ exceeds every threshold, only the equilibria of zero energy remain, and the degree counts extrema minus saddles on the nodal set of $W$. In this sense the family $q^2+1+\varepsilon\nabla W(\Ima q)$ filters the equilibria by energy, from the most pronounced to the most neutral.

\subsection{Three classical energies}\label{subsec:anisotropy-examples}
We now consider the biaxial, octupolar and cubic energies. Each of them realizes the equality case of~\eqref{eq:anisotropy-morse} for its degree, so the extremal fields of Theorem~\ref{thm:sphere-bounded-count} for $D\leq3$ need no appeal to the general existence result of Kozhasov.

For $D=1$, the biaxial energy $W(x)=\frac12(k_1x_1^2+k_2x_2^2+k_3x_3^2)$ with $k_1<k_2<k_3$ has $h(x)=Mx$ with $M=\operatorname{diag}(k_1,k_2,k_3)$. Its six equilibrium directions are $\pm e_1$, the easy axis, which gives minima, $\pm e_2$, the intermediate axis, which gives saddles, and $\pm e_3$, the hard axis, which gives maxima. This is the configuration of Proposition~\ref{prop:sphere-linear}, whose roots lie over these six points with local degrees $+1,-1,+1$ in the order of the eigenvalues. The uniaxial energy, with two equal coefficients, has a circle of equilibria and falls outside the dictionary, in accordance with the non-degeneracy hypothesis.

For $D=2$, the octupolar energy $W_o(x)=x_1x_2x_3$ has gradient $h=(x_2x_3,x_1x_3,x_1x_2)$, the field of Theorem~\ref{thm:sphere-xyz}. Energies of this kind appear in the description of octupolar order in liquid crystals \cite{gaetavirga}. The proof of Theorem~\ref{thm:sphere-xyz} shows that its equilibrium directions are the six saddles $\pm e_a$, of energy zero, and the eight points $s/\sqrt3$ with $s\in\{\pm1\}^3$, of energy $s_1s_2s_3/(3\sqrt3)$. The four with $s_1s_2s_3=1$ are maxima, and the four others are minima. Thus $14=2(2^2+2+1)$ equilibria occur. The threshold $2/(3|\varepsilon|)$ of Proposition~\ref{prop:anisotropy-roots} crosses the energy $1/(3\sqrt3)$ of the extrema at $|\varepsilon|=2\sqrt3$, and the passage of the Brouwer degree from $2$ to $-6$ in Theorem~\ref{thm:sphere-xyz} is the departure of the eight extrema, which leaves the six saddles on the nodal set.

For $D=3$, the leading term of the cubic magnetocrystalline energy is $K_1W_c$ with
\begin{equation}\label{eq:cubic-energy}
W_c(x)=x_1^2x_2^2+x_2^2x_3^2+x_3^2x_1^2
\end{equation}
where $K_1$ is the first cubic anisotropy constant \cite{chikazumi}. Its sign decides between easy axes along the cube edges, as for iron, and along the cube diagonals, as for nickel.

\begin{proposition}[Cubic anisotropy]\label{prop:cubic-anisotropy}
The energy $E_c=W_c|_{\mathbb S^2}$ is a Morse function with exactly $26=2(3^2+3+1)$ critical points. They are
\begin{enumerate}
\item the six points $\pm e_a$, of energy $0$, which are minima,
\item the twelve points $(\pm e_a\pm e_b)/\sqrt2$ with $a<b$, of energy $\frac14$, which are saddles,
\item the eight points $s/\sqrt3$ with $s\in\{\pm1\}^3$, of energy $\frac13$, which are maxima.
\end{enumerate}
In particular $e(\nabla W_c)=13$, and $\nabla W_c$ realizes the equality case of Theorem~\ref{thm:sphere-bounded-count} for $D=3$.
\end{proposition}
\begin{proof}
The gradient is $h_c(x)=2\bigl(x_1(x_2^2+x_3^2),x_2(x_1^2+x_3^2),x_3(x_1^2+x_2^2)\bigr)$. For $J\in\mathbb S^2$ one has $J_2^2+J_3^2=1-J_1^2$ and its analogues, so the eigenvector equation $h_c(J)=\mu J$ becomes $2J_a(1-J_a^2)=\mu J_a$ for $a=1,2,3$. Every non-zero coordinate therefore has the same square $1-\mu/2$. When $k\in\{1,2,3\}$ coordinates are non-zero, this square equals $1/k$, and the three values of $k$ give the three families of the statement, with $\mu=2(1-1/k)$ and energy $\mu/4$ by Lemma~\ref{lem:anisotropy}. This accounts for all critical points.
The signed permutations of the coordinates preserve $W_c$ and act transitively on each family, so it suffices to compute $K_J$ at one point of each. The derivative $Dh_c(J)$ has diagonal entries $2(J_b^2+J_c^2)$ and off-diagonal entries $4J_aJ_b$, where $\{a,b,c\}=\{1,2,3\}$. At $J=e_1$ one has $\mu=0$ and $K_J=\operatorname{diag}(2,2)$ on the span of $e_2,e_3$, so $J$ is a minimum. At $J=(e_1+e_2)/\sqrt2$ one has $\mu=1$, and the tangent vectors $\tau_1=(e_1-e_2)/\sqrt2$ and $\tau_2=e_3$ satisfy $Dh_c(J)\tau_1=-\tau_1$ and $Dh_c(J)\tau_2=2\tau_2$. Hence $K_J=\operatorname{diag}(-2,1)$, and $J$ is a saddle. At $J=(1,1,1)/\sqrt3$ one has $\mu=\frac43$, and every entry of $Dh_c(J)$ equals $\frac43$. This matrix annihilates the tangent plane, so $K_J=-\frac43\mathrm{Id}$ and $J$ is a maximum. The count $6+8-12=2$ agrees with~\eqref{eq:anisotropy-morse}.
\end{proof}

With $K_1>0$, the minima of $K_1E_c$ are the cube edges, the saddles are the face diagonals and the maxima are the cube diagonals. With $K_1<0$ minima and maxima exchange their roles, while the saddles remain. The quaternionic statements below are insensitive to this sign, since the local degree of an extremum does not distinguish a minimum from a maximum.

\begin{corollary}[Explicit extremal roots from cubic anisotropy]\label{cor:cubic-anisotropy-roots}
For every sufficiently small $\varepsilon\neq0$, the map $q^4+1+\varepsilon\nabla W_c(\Ima q)$ has exactly $52$ roots, all simple. The $28$ roots over the cube edges and cube diagonals have local degree $+1$, and the $24$ roots over the face diagonals have local degree $-1$.
\end{corollary}
\begin{proof}
This is Proposition~\ref{prop:anisotropy-roots} with $m=2$, together with Proposition~\ref{prop:cubic-anisotropy}.
\end{proof}

The stacked quartic~\eqref{eq:en-quartic-example} has the $132$ simple roots and coefficient stability given in Corollary~\ref{cor:en-quartic}.

\subsection{Reach of the dictionary}\label{subsec:anisotropy-reach}
The dictionary transfers counts, tangent indices and energy thresholds. The quaternionic variable carries no physical meaning of its own, and dynamics lie outside the correspondence. Magnetic anisotropy energies are invariant under $J\mapsto-J$, by time-reversal symmetry, so their degree $D+1$ is even and $D$ is odd. Lemma~\ref{lem:eigen-antipodal} then gives antipodal equilibria the same energy, and the count $2e(h)$ is congruent to two modulo four. Octupolar energies have odd degree, and antipodal equilibria carry opposite energies, which explains the four maxima and four minima of $W_o$.

For $D\geq3$, the family $q^2+1+\varepsilon\nabla W(\Ima q)$ also has explicit energy thresholds. Assume that $E$ is a Morse function and $\varepsilon\neq0$. By~\eqref{eq:eigen-quadratic-planar}, the roots over an equilibrium direction $J$ are $q=-\frac12\varepsilon\mu_Js^{(D-1)/2}+\sqrt{s}\,J$, where $\mu_J=(D+1)E(J)$ and
\begin{equation}\label{eq:anisotropy-higher-radius}
s>0,\qquad \frac{\varepsilon^2\mu_J^2}{4}s^{D-1}-s+1=0.
\end{equation}
Since $\mathcal E(h)$ is finite and each radial polynomial is non-zero, every fixed $\varepsilon\neq0$ has finitely many roots. An isolating ball therefore exists without a global properness assumption. For $E(J)=0$ there is one root, $q=J$. For $E(J)\neq0$, there are two roots when $0<|\varepsilon|<\varepsilon_J^*$, one double radial solution at equality, and no roots above the threshold
\begin{equation}\label{eq:anisotropy-higher-threshold}
\varepsilon_J^*=\frac{2}{(D+1)|E(J)|}\left(\frac{(D-2)^{D-2}}{(D-1)^{D-1}}\right)^{1/2}.
\end{equation}
Indeed, writing $c=\varepsilon^2\mu_J^2/4$, the polynomial $cs^{D-1}-s+1$ is convex on $s>0$ and has a unique minimum there. This minimum vanishes exactly when $s=(D-1)/(D-2)$ and $c=(D-2)^{D-2}/(D-1)^{D-1}$. At a root, the planar determinant~\eqref{eq:eigen-quadratic-determinant} equals $4s-4(D-1)cs^{D-1}=4\bigl((D-1)-(D-2)s\bigr)$. Lemma~\ref{lem:eigen-reduction} therefore gives opposite local degrees for the two roots, the root of smaller radius having local degree $\operatorname{ind}_J(V_h)$. Away from the thresholds, the signed sum of all roots is $\sum_{J\in\mathcal E(h),\,E(J)=0}\operatorname{ind}_J(V_h)$ for every $\varepsilon\neq0$. For each fixed $\varepsilon$ away from the thresholds, this sum is the Brouwer degree on any such isolating ball. The unsigned count changes at the thresholds, while this signed sum remains constant.

At the threshold the two roots merge at an ordinary fold. The tangent block of Lemma~\ref{lem:eigen-reduction} is invertible, and the full planar Jacobian determinant is
\begin{equation}\label{eq:anisotropy-fold-jacobian}
\Delta_J(x,y)=4x^2+4y^2+2D\varepsilon\mu_Jxy^{D-1}.
\end{equation}
At the critical root, the planar derivative has kernel spanned by $(y,x)$. Using $\varepsilon\mu_Jy^{D-1}=-2x$ and $(D-1)x^2=y^2$, the derivative of $\Delta_J$ along this kernel equals $-8(D-2)xy\neq0$. Along $\CC_J$ near the critical root, the lower-left block of $DG$ vanishes, and its tangent block equals $\varepsilon y^{D-1}K_J$ at the root. This block stays invertible on the slice near the root, so $\det DG$ there is $\Delta_J$ times a non-vanishing factor. The kernel of $DG$ lies in this slice. The derivative of $\det DG$ along the planar kernel at the critical root is therefore a non-zero multiple of $-8(D-2)xy$. The corank-one fold criterion then applies to the four-dimensional map \cite[Chapter~III, \S4]{golubitskyguillemin}.

As $\varepsilon$ tends to zero, the root of smaller radius approaches $J$ and the root of larger radius escapes to infinity. The pairs over the equilibria of non-zero energy thus disappear in decreasing order of $|E(J)|$ as the coupling grows. For the cubic energy~\eqref{eq:cubic-energy}, the thresholds are $\frac34$ over the cube diagonals and $1$ over the face diagonals. For $D=2$, the same scalar equation becomes $(1-c)s=1$, so the threshold is $c=1$, that is, $|\varepsilon\mu_J|=2$ as in Proposition~\ref{prop:quad-wall}. The radius then tends to infinity as $c$ increases to one, while the collision radius for $D\geq3$ is finite.

The second root over each equilibrium of non-zero energy lies far from the unit sphere when $\varepsilon$ is small. It is therefore outside the tube of Theorem~\ref{thm:sphere-splitting}, and it does not occur in Theorem~\ref{thm:en-splitting}, where the perturbation has degree less than that of the leading power.

\section{An outlook on orientation finding in diffusion MRI}\label{sec:mri-orientation}
The equilibrium problem on $\mathbb S^2$ also appears in diffusion magnetic resonance imaging. A fitted fiber orientation distribution describes the directional organization of tissue within a voxel, and its local maxima provide candidate fiber directions for tractography \cite{tournier2007}. An antipodally symmetric spherical-harmonic expansion gives a finite representation of such a profile. Polynomial extraction and classification of spherical critical points have been developed for diffusion MRI by Ghosh, Tsigaridas, Mourrain and Deriche \cite{ghosh2013}. The eigenvector reduction of Section~\ref{sec:eigen} connects this problem with spherical root splitting, and its count and index formulas provide checks on an enumeration of candidate directions.

Let $F$ be a real antipodally symmetric profile whose highest non-zero spherical-harmonic order is an even integer $L\geq2$. Write $F(J)=\sum_{\ell=0,2,\ldots,L}\mathcal H_\ell(J)$, where $\mathcal H_\ell$ is a real homogeneous harmonic polynomial of degree $\ell$. The homogeneous polynomial
\begin{equation}\label{eq:mri-homogeneous-lift}
W_F(v)=\sum_{\ell=0,2,\ldots,L}|v|^{L-\ell}\mathcal H_\ell(v)
\end{equation}
has degree $L$ and restricts to $F$ on $\mathbb S^2$. The exponent $(L-\ell)/2$ is a non-negative integer, and $|v|^{L-\ell}=(v_1^2+v_2^2+v_3^2)^{(L-\ell)/2}$, so the lift is polynomial in the three real coordinates. Lemma~\ref{lem:anisotropy} identifies the critical directions of $F$ with the unit eigenvectors of $\nabla W_F$, whose degree is $D=L-1$. They are therefore the real solutions of
\begin{equation}\label{eq:mri-critical-system}
\nabla W_F(J)=\mu J,\qquad |J|^2=1,\qquad J\in\RR^3,\quad\mu\in\RR.
\end{equation}
At a non-degenerate solution, two positive eigenvalues of $K_J$ give a minimum, two negative ones give a maximum, and opposite signs give a saddle. Antipodal maxima represent a single unoriented fiber direction.

For a Morse profile $F$, Theorem~\ref{thm:sphere-bounded-count} bounds the total number of critical points by $2(L^2-L+1)$, while~\eqref{eq:anisotropy-morse} gives $n_{\min}+n_{\max}-n_{\mathrm{sad}}=2$. Adding these relations gives $n_{\min}+n_{\max}\leq L^2-L+2$. Antipodal minima give $n_{\min}\geq2$, hence $n_{\max}\leq L(L-1)$ and at most $L(L-1)/2$ unoriented peak directions. Equality in the peak bound requires $n_{\min}=2$ and exactly $2(L^2-L+1)$ critical points. The simultaneous realization of these conditions is not proved here. At order eight there are at most $114$ critical points, or $57$ antipodal pairs, and at most $28$ unoriented peak directions. This ceiling concerns the fitted polynomial profile. Its coefficients depend on the response model and, in constrained spherical deconvolution, on the non-negativity constraint \cite{tournier2007}.

Let $\gamma_J>0$ be the smallest absolute eigenvalue of the tangent Hessian $K_J$. For a coefficient variation $\delta F$ and the corresponding first variation $\delta J\in T_J\mathbb S^2$ of a non-degenerate critical direction, differentiating the critical equation gives
\begin{equation}\label{eq:mri-critical-sensitivity}
\delta J=-K_J^{-1}\nabla_{\mathbb S^2}\delta F(J),\qquad |\delta J|\leq\frac{|\nabla_{\mathbb S^2}\delta F(J)|}{\gamma_J},
\end{equation}
where $\nabla_{\mathbb S^2}$ denotes the spherical gradient. For fixed $L$, the algebraic ceiling is the same for every Morse profile. Small $\gamma_J$ increase location sensitivity and make numerical verification of non-degeneracy more demanding.

A rigorously verified list of $2(L^2-L+1)$ distinct critical points of $W_F|_{\mathbb S^2}$ is complete under the Morse hypothesis, since the homogeneous lift~\eqref{eq:mri-homogeneous-lift} has degree $L$. Disjoint enclosures with verified existence and uniqueness of a root, together with verified Hessian eigenvalue signs, provide the local checks. During numerical critical-point tracking, the index sum supplies a necessary consistency check at smaller counts. An extremum and a saddle omitted together leave the sum unchanged, and repeated searches with the same result supply numerical evidence. A certificate at a smaller count requires an independent exhaustion of~\eqref{eq:mri-critical-system} over the reals. Peak-only searches leave the saddles outside their enumeration. Searching for all zeros of $V_{\nabla W_F}$ and classifying them with the projected Hessian $K_J$ of~\eqref{eq:eigen-derivative} permits this fuller check.

Following maxima together with nearby saddles and monitoring $\gamma_J$ offers a way to study changes in peak structure as the fitted coefficients vary. These checks concern the fitted profile and its critical directions. Their value for fiber reconstruction requires assessment on acquired data, with the response model, angular resolution and measurement noise taken into account.

\section*{Methods and AI statement}\label{sec:ai-usage}
The paper began with the author's logical and categorical investigation into proofs of the Fundamental Theorem of Algebra and the boundary between proofs that are available in classical logic and those that are available in an intuitionistic framework. 

Frontier large language models, namely OpenAI's GPT-6-Astra and GPT-6.1-Sol, and Kimi-K3, were then prompted with this categorical framework to search for results at the boundary where classical formulation may obscure results that are available constructively. The categorical and logical apparatus is omitted from the present paper, which focuses entirely on the resulting theorems.

Specific mathematical contributions due to GPT-6.1 Sol include the uniform sign estimates and planar Jacobian calculation in Corollary~\ref{cor:en-quartic}, and the critical-direction sensitivity estimate in Section~\ref{sec:mri-orientation}.

The author assumes full responsibility for the integrity, citations, and content of the entire work.

\bibliographystyle{smfalpha}
\bibliography{quaternionic}
\end{document}